\documentclass[a4paper,12pt]{article}

\usepackage{amscd}
\usepackage{amsfonts}
\usepackage{amsmath}
\usepackage{amssymb}
\usepackage{amsthm}
\usepackage{ascmac}
\usepackage{here}
\usepackage{mathrsfs}
\usepackage{slashbox}
\usepackage{txfonts}
\usepackage[all]{xy}

\allowdisplaybreaks

\newcommand{\C}{\mathbb{C}}
\newcommand{\F}{\mathbb{F}}
\newcommand{\N}{\mathbb{N}}
\newcommand{\Q}{\mathbb{Q}}
\newcommand{\R}{\mathbb{R}}
\newcommand{\Z}{\mathbb{Z}}

\newcommand{\Aff}{\mathbb{A}}

\newcommand{\A}{\mathscr{A}}
\newcommand{\B}{\mathscr{B}}

\newcommand{\Fil}{\mathscr{F}}

\newcommand{\Leb}{\mathscr{L}}

\newcommand{\Nor}{\mathscr{N}}
\newcommand{\Opn}{\mathscr{O}}

\newcommand{\U}{\mathscr{U}}

\theoremstyle{plain}
\newtheorem{thm}{Theorem}[section]
\newtheorem{lmm}[thm]{Lemma}
\newtheorem{prp}[thm]{Proposition}
\newtheorem{crl}[thm]{Corollary}
\theoremstyle{definition}
\newtheorem{dfn}[thm]{Definition}
\newtheorem{rmk}[thm]{Remark}
\newtheorem{exm}[thm]{Example}

\def\ens#1{\{ #1 \}}
\def\Ens#1{\left\{ #1 \right\}}
\def\set#1#2{\{ #1 \mid #2 \}}
\def\Set#1#2{\left\{ #1 \ \middle| \ #2 \right\}}
\def\t#1{\text{\rm #1}}

\def\v#1{| #1 |}
\def\vv#1{\left| #1 \right|}
\def\n#1{\| #1 \|}
\def\nn#1{\left\| #1 \right\|}

\title{Spectral Theory for $p$-Adic Operators}
\author{Tomoki Mihara\footnote{E-mail address: mihara@ms.u-tokyo.ac.jp}}
\date{}

\begin{document}

\maketitle
\begin{abstract}
We establish an algorithm for a criterion of the diagonalisability of a matrix over a local field by a unitary matrix. For this sake, we define the notion of normality of a $p$-adic operator, and give several criteria for the normality. We study the relation between the normality and the reduction. In the finite dimensional case, the normality of an operator is equivalent to the diagonalisability of a matrix by a unitary matrix. Therefore we also study the relation between the diagonalisability and the reduction. For example, we show that the diagonalisation of the reduction gives a partition of unity corresponding to the reduction of the spectrum, which gives a functorial lift of the eigenspace decomposition of the reduction.
\end{abstract}
\tableofcontents

\section{Introduction}
\label{Introduction}

Let $k$ be a complete valuation field. The main aim of this paper is to establish an algorithm for a criterion of the diagonalisability of a matrix over $k$ by a unitary matrix. For this sake, we define the normality of a $k$-linear bounded operator in Definition \ref{normal} and establish several criteria for the normality using information of the reduction in Corollary \ref{n implies rn 3}, Corollary \ref{hfc for a reductively infinite operator 2}, and Corollary \ref{reduction and functional calculus 2}. A $k$-linear bounded operator is said to be {\it normal} if and only if it admits the continuous functional calculus (Definition \ref{rcfc}), and the notion of the normality of a $k$-linear bounded operator on a finite dimensional $k$-vector space is equivalent to that of the diagonalisability of a matrix over $k$ by a unitary matrix (Proposition \ref{n <-> ud}). Our formulation of the normality is a generalisation of the other formulation of the normality in \cite{Koc13} \S 2. A straightforward method for determining the normality includes pure $p$-adic analysis. In order to avoid from the analysis, we study a relation between reductions, spectra, and functional calculi. Since a functional calculus of an operator over the residue field is just the diagonalisation of a matrix, a reduction helps us to reduce the analysis to pure ring theory. In the case where $k$ is a local field, we give an explicit algorithm for a criterion of the diagonalisability of a matrix over $k$ by a unitary matrix including a ``repetition of reductions'' in Theorem \ref{unitarily diagonalisable}. In particular, it yields a sufficient condition for the diagonalisability of a matrix.

\vspace{0.1in}
A {\it functional calculus} is a substitution of a bounded operator to a function on its spectrum. For example, the holomorphic functional calculus is a substitution to a rigid analytic function, and the continuous functional calculus is a substitution to a continuous function. They are the extensions of the substitution to a polynomial function. In order to justify substitutions to such functions, we need a formulation of a spectral decomposition, which is a generalised notion of the eigenspace decomposition of a diagonalisable matrix. In order to deal with a continuous functional analysis, we establish a theory of a spectral decomposition. A formulation of a spectral decomposition of an analytic bounded operator with compact spectrum has been already given by Vishik in his paper \cite{Vis85}, but it just yields kind of the locally analytic functional calculus analogous to Riesz functional calculus. Under several suitable conditions for a bounded operator, there are direct extensions of the result of Vishik by Dodzi Attimu and Toka Diagana in \cite{AD09} and by Richard Lance Baker in \cite{Bak12}, but the direction of their works is completely different from that of ours. We would like to remove the assumptions that a bounded operator is analytic and that its spectrum is compact, which are less useful for a continuous functional calculus. Such a non-trivial generalisation can not be obtained if one extends Vishik's spectral theory in a direct way. Indeed, the locally analytic functional calculus in his work uses Shnirel'man integral, which is a non-Archimedean analogue of Cauchy integral of an analytic function introduced in \cite{Shn38}. The assumptions that a bounded operator is analytic and that its spectrum is compact are essential for the use of Shnirel'man integral. We briefly explain how we establish a spectral decomposition of a bounded operator. Let $A$ be a bounded operator on a Banach $k$-vector space. In the case where $k$ is a local field, the spectrum $\sigma$ of $A$ is automatically compact. In this case, we formulate the measurability of a clopen subset of $\sigma$ in Definition \ref{measurable}, and establish a method to give a spectral decomposition as an integration along $\sigma$ without the assumption that $A$ is analytic in Theorem \ref{naivety}. Even if $k$ is not a local field, we also give a sufficient condition of the integrability along the spectrum using affinoid algebras without the assumption of the compactness of $\sigma$ in Theorem \ref{hfc for a reductively infinite operator}.

\section{Preliminaries}
\label{Preliminaries}

We recall basic notions and conventions in Banach spaces and Banach algebras over a complete valuation field in \S \ref{Banach Spaces and Banach Algebras}. We follow the conventions in \cite{BGR84}, except that we deal also with non-commutative Banach algebras. As important examples of Banach algebras, we introduce rings of ``rigid continuous functions'' and operator algebras in \S \ref{Rigidity of Continuous Functions} and \ref{Operator Algebra} respectively.

\subsection{Banach Spaces and Banach Algebras}
\label{Banach Spaces and Banach Algebras}

Throughout this paper, we assume that rings and algebras are unital and associative. We do not assume the commutativity or the inequality $1 \neq 0$. A field is assumed to be a commutative ring with $0 \neq 1$. In particular, the zero ring $O$ is excluded from the definition of a field. A {\it valuation field} is a field $k$ endowed with a map $\v{\cdot} \colon k \to [0, \infty)$ satisfying the following:
\begin{itemize}
\item[(i)] The inequality $\v{c-c'} \leq \max \ens{\v{c}, \v{c'}}$ holds for any $(c,c') \in k \times k$.
\item[(ii)] The equality $\v{cc'} = \v{c} \ \v{c'}$ holds for any $(c,c') \in k \times k$.
\item[(iii)] The equality $\v{c} = 0$ holds for a $c \in k$ if and only if $c = 0$.
\item[(iv)] The equality $\v{1} = 1$ holds.
\end{itemize}
We always endow a valuation field $k$ with the topology given by the ultrametric $k \times k \to [0,\infty) \colon (c,c') \mapsto \v{c-c'}$. A valuation field is said to be {\it complete} if its underlying ultrametric is complete, and is said to be {\it discrete} if the image $\v{k^{\times}}$ of $k^{\times}$ by $\v{\cdot}$ is a subgroup of the multiplicative group $(0,\infty)$ isomorphic to $\Z$. Let $k$ be a valuation field. We denote by $O_k \subset k$ the clopen subring consisting of elements $c$ with $\v{c} \leq 1$, and call it {\it the valuation ring of $k$}. The valuation ring $O_k$ is a local ring. We denote by $\overline{k}$ the discrete field obtained as the quotient of $O_k$ by its unique maximal ideal $m_k \subset O_k$, and call it {\it the residue field of $k$}. We regard $\overline{k}$ a complete valuation field with respect to the trivial valuation. A valuation field is said to be a {\it local field} if it is a complete discrete valuation field whose residue field is a finite field.

\begin{exm}
Every field $\F$ admits a valuation $\v{\cdot} \colon \F \to [0,\infty)$ called {\it the trivial valuation} given by setting $\v{0} \coloneqq 0$ and $\v{c} \coloneqq 1$ for any $c \in \F^{\times}$.
\end{exm}

Let $k$ be a complete valuation field. For a $k$-vector space $W$, a {\it seminorm} on $W$ is a map $\n{\cdot} \colon W \to [0,\infty)$ satisfying the following:
\begin{itemize}
\item[(i)] The inequality $\n{v-v'} \leq \max \ens{\n{v}, \n{v'}}$ holds for any $(v,v') \in W \times W$.
\item[(ii-a)] The equality $\n{cv} = \v{c} \ \n{v}$ holds for any $(c,v) \in k \times W$.
\end{itemize}
For a $k$-vector space $W$, a {\it norm on $W$} is a seminorm on $W$ satisfying the following:
\begin{itemize}
\item[(iii)] The equality $\n{v} = 0$ holds for a $v \in W$ if and only if $v = 0$.
\end{itemize}
A {\it normed $k$-vector space} is a $k$-vector space $V$ endowed with a norm. We always endow a normed $k$-vector space $V$ with the topology given by the ultrametric $V \times V \to [0,\infty) \colon (v,v') \mapsto \n{v-v'}$. A normed $k$-vector space is said to be a {\it Banach $k$-vector space} if its underlying ultrametric is complete. In particular, $k$ itself is a Banach $k$-vector space.

\begin{exm}
Suppose that the valuation of $k$ is trivial. Every $k$-vector space $W$ admits a norm $\n{\cdot} \colon W \to [0,\infty)$ called {\it the trivial norm} given by setting $\n{0} = 0$ and $\n{v} = 1$ for any $v \in W$. A normed $k$-vector space is said to be a {\it discrete $k$-vector space} if it is endowed with the trivial norm. Every discrete $k$-vector space is a Banach $k$-vector space, while a Banach $k$-vector space is not necessarily a discrete $k$-vector space. For example, $k^{\N}$ is a Banach $k$-algebra with respect to the norm $\n{\cdot} \colon k^{\N} \to [0,\infty)$ given by setting $\n{(0)_{n \in \N}} \coloneqq 0$ and $\n{(a_n)_{n \in \N}} \coloneqq \exp(- \sup \set{n \in \N}{a_{n_0} = 0, {}^{\forall} n_0 \in \N \cap [0,n)})$ for each $(a_n)_{n \in \N} \in k^{\N} \backslash \ens{(0)_{n \in \N}}$, and is not a discrete $k$-vector space.
\end{exm}

Let $V$ be a normed $k$-vector space. For each $r > 0$, we denote by $V(r) \subset V$ (resp.\ $V(r-) \subset V$) the clopen $O_k$-submodule consisting of elements $v$ with $\n{v} \leq r$ (resp.\ $\n{v} < r$), and call it {\it the closed ball of $V$ of radius $r$} (resp.\ the open ball of $V$ radius $r$). The quotient $\overline{V} \coloneqq V(1)/V(1-)$ is a normed $\overline{k}$-vector space with respect to the trivial norm. We call it {\it the reduction of $V$}. For each $v \in V(1)$, we denote by $v + V(1-) \in \overline{V}$ the image of $v$, which is the coset $\set{v+v'}{v' \in V(1-)} \subset V(1)$ by the set-theoretical definition of the quotient set.

\begin{exm}
Let $X$ be a topological space. A continuous function $f \colon X \to k$ is said to {\it vanish at infinity} if for any $\epsilon > 0$, there is a compact subspace $K \subset X$ with $\v{f(x)} < \epsilon$ for any $x \in X \backslash K$. The $k$-vector space $\t{C}_0(X,k)$ of continuous functions $X \to k$ vanishing at infinity forms a Banach $k$-vector space with respect to the supremum norm $\n{\cdot} \colon \t{C}_0(X,k) \to [0,\infty)$ given by setting $\n{f} \coloneqq \sup_{x \in X} \v{f(x)}$ for each $x \in X$. If the valuation of $k$ is trivial, then $\t{C}_0(X,k)$ is a discrete $k$-vector space. In addition if $X$ is a discrete topological space $I$, then the canonical embedding $k^{\oplus I} \hookrightarrow k^I$ gives an isometric $k$-linear isomorphism $k^{\oplus I} \to \t{C}_0(I,k)$ with respect to the trivial norm on $k^{\oplus I}$, and hence we identify $\t{C}_0(I,k)$ and $k^{\oplus I}$.
\end{exm}

For topological spaces $X$ and $Y$, we denote by $\t{C}(X,Y)$ the $k$-vector space of continuous maps $X \to Y$. Let $X$ be a compact topological space. Every continuous function $X \to k$ vanishes at infinity by definition, and hence $\t{C}(X,k)$ coincides with the underlying $k$-vector space of $\t{C}_0(X,k)$. Identifying $\t{C}(X,k)$ with $\t{C}_0(X,k)$, we regard $\t{C}(X,k)$ as a Banach $k$-vector space. In the following, we regard every set $I$ as a discrete topological space so that the convention $\t{C}_0(I,k)$ makes sense.

\begin{prp}
\label{reduction & direct product}
For any set $I$, there is a natural $\overline{k}$-linear isomorphism $\overline{\t{C}_0(I,k)} \to \overline{k}{}^{\oplus I}$.
\end{prp}

\begin{proof}
The continuous $O_k$-linear homomorphism $\t{C}_0(I,k) \to \overline{k}{}^I \colon f \mapsto (f(i) + m_k)_{i \in I}$ induces an injective $\overline{k}$-linear homomorphism $\overline{\t{C}_0(I,k)} \to \overline{k}{}^I$, whose image coincides with $\overline{k}{}^{\oplus I}$ by the definition of $\t{C}_0(I,k)$.
\end{proof}

Let $V_1$ and $V_2$ be normed $k$-vector spaces. A $k$-linear homomorphism $f \colon V_1 \to V_2$ is said to be {\it bounded} if $\set{\n{f(v)}}{v \in V(1)}$ is a bounded subset of $[0,\infty)$, and is said to be {\it submetric} if the inequality $\n{f(v)} \leq \n{v}$ holds for any $v \in V_1$. The $k$-vector space $\t{Hom}_k^{\t{cont}}(V_1,V_2)$ of bounded $k$-linear homomorphisms forms a normed $k$-vector space with respect to {\it the operator norm} $\n{\cdot} \colon \t{Hom}_k^{\t{cont}}(V_1,V_2) \to [0,\infty)$ given by setting $\n{f} \coloneqq \sup_{v \in V_1(1)} \n{f(v)}$ for each $f \in \t{Hom}_k^{\t{cont}}(V_1,V_2)$. If $V_2$ is a Banach $k$-vector space, then so is $\t{Hom}_k^{\t{cont}}(V_1,V_2)$.

\begin{dfn}
Let $V$ be a Banach $k$-vector space. A discrete subset $I \subset V$ is said to be an {\it orthonormal Schauder basis} of $V$ if the $k$-linear homomorphism $k^{\oplus I} \to V$ given by the inclusion $I \hookrightarrow V$ extends to an isometric isomorphism $\t{C}_0(I,k) \to V$ through the natural injective $k$-linear homomorphism $k^{\oplus I} \subset \t{C}_0(I,k)$. If $V$ is of finite dimension, then we simply call $I$ an {\it orthonormal basis}.
\end{dfn}

An orthonormal basis of a finite-dimensional Banach $k$-vector space is a $k$-linear basis of the underlying $k$-vector space, but the converse does not hold in general. There are criteria for the orthonormality in \cite{BGR84} 2.5.1/3 and \cite{BGR84} 2.5.2/2. We note that the original notion of being strictly Cartesian is defined for a normed $k$-vector space. For example, see \cite{BGR84} 2.5.2/1.

\begin{dfn}
A Banach $k$-vector space is said to be {\it strictly Cartesian} if it admits an orthonormal Schauder basis, or equivalently if it is isometrically isomorphic to the Banach $k$-vector space $\t{C}_0(I,k)$ for some discrete topological space $I$.
\end{dfn}

\begin{exm}
Let $n \in \N$. Then $k^n$ is a finite-dimensional strictly Cartesian Banach $k$-vector space with respect to the norm given by the bijective $k$-linear homomorphism $\t{C}_0(\N \cap [1,n],k) \to k^n$ associated to the canonical basis. If one replaces $k^n$ by another equivalent norm, then $k^n$ is not necessarily strictly Cartesian. Every subspace of a finite-dimensional strictly Cartesian Banach space is again strictly Cartesian by \cite{BGR84} 2.5.1/4.
\end{exm}

\begin{exm}
Suppose that the valuation of $k$ is trivial. A Banach $k$-vector space is strictly Cartesian if and only if it is a discrete $k$-vector space, and every $k$-linear basis of a discrete $k$-vector space is an orthonormal Schauder basis.
\end{exm}

\begin{dfn}
A Banach $k$-vector space is said to be {\it of countable type} if it contains dense $k$-vector subspace of countable dimension.
\end{dfn}

If $k$ is separable, then the notion of a Banach $k$-vector space of countable type is equivalent to that of a separable Banach $k$-vector space. It is remarkable that a Banach space of countable type always admits an equivalent norm with respect to which it admits an orthonormal Schauder basis by \cite{BGR84} 2.7.2/8. Such a countability can be removed if the valuation of $k$ is discrete. For more details, see \cite{Mon70} IV.3, \cite{BGR84} 2.5.3/11 and 2.4.4/2, or \cite{FP04} 1.2.2/(1).

\begin{exm}
The $\Q_p$-Banach spaces $\t{C}(\Z_p,\Q_p)$, $\t{C}_0(\Q_p,\Q_p)$, and $\C_p$ are strictly Cartesian $\Q_p$-Banach spaces of countable type.
\end{exm}

\begin{exm}
Suppose that the valuation of $k$ is trivial. A discrete $k$-vector space is of countable type if and only if its underlying $k$-vector space is of countable dimension.
\end{exm}

For a $k$-algebra $A$, a {\it seminorm on $A$} (resp.\ a {\it norm on $A$}) is a seminorm (resp.\ a norm) $\n{\cdot}$ on the underlying $k$-vector space of $A$ satisfying the following:
\begin{itemize}
\item[(ii-b)] The inequality $\n{ff'} \leq \n{f} \ \n{f'}$ for any $(f,f') \in A \times A$.
\item[(iv)] The equality $\n{1} = 1$ holds.
\end{itemize}
For a $k$-algebra $A$, a seminorm $\n{\cdot}$ on $A$ is said to be {\it power-multiplicative} if it satisfying the following:
\begin{itemize}
\item[(ii-c)] The equality $\n{f^n} = \n{f}^n$ holds for any $(f,n) \in A \times \N$.
\end{itemize}
For a $k$-algebra $A$, a seminorm $\n{\cdot}$ on $A$ is said to be {\it multiplicative} if it satisfying the following:
\begin{itemize}
\item[(ii)] The equality $\n{ff'} = \n{f} \ \n{f'}$ holds for any $(f,f') \in A \times A$.
\end{itemize}
A {\it normed $k$-algebra} is a $k$-algebra $\A$ endowed with a norm. A normed $k$-algebra is said to be a {\it Banach $k$-algebra} if its underlying normed $k$-vector space is a Banach $k$-vector space. A Banach $k$-algebra is said to be {\it uniform} if its norm is power-multiplicative. Let $\A$ be a normed $k$-algebra. We denote by $\A(1), \A(1-) \subset \A$ the closed unit ball and the open unit ball of the underlying normed $k$-vector space of $\A$. Then $\A(1)$ is a clopen $O_k$-subalgebra of $O_k$, and $\A(1-)$ is a clopen ideal of $\A(1)$. Therefore the reduction $\overline{\A} \coloneqq \A(1)/\A(1-)$ of the underlying normed $k$-vector space of $\A$ admits a natural structure of a discrete $\overline{k}$-algebra. We call $\overline{\A}$ {\it the reduction of $\A$}.

\begin{exm}
Let $X$ be a topological space. A continuous function $f \colon X \to k$ is said to be {\it bounded} if $\set{\v{f(x)}}{x \in X}$ is a bounded subset of $[0,\infty)$. The $k$-algebra $\t{C}_{\t{bd}}(X,k)$ of bounded continuous functions $X \to k$ forms a uniform Banach $k$-algebra with respect to the power-multiplicative $\n{\cdot} \colon \t{C}_{\t{bd}}(X,k) \to [0,\infty)$ called {\it the supremum norm} given by setting $\n{f} \coloneqq \sup_{x \in X} \v{f(x)}$ for each $x \in X$. If the valuation of $k$ is trivial and $X$ is a discrete topological space $I$, then $\t{C}_{\t{bd}}(I,k)$ coincides with $k^I$ endowed with the trivial norm.
\end{exm}

Let $X$ be a topological space $X$. Suppose that $X$ is compact or the valuation of $k$ is trivial. Then every continuous function $X \to k$ is bounded by the maximum modulus principle, and hence $\t{C}(X,k)$ coincides with the underlying Banach $k$-vector space of $\t{C}_{\t{bd}}(X,k)$. Identifying $\t{C}(X,k)$ and $\t{C}_{\t{bd}}(X,k)$, we regard $\t{C}(X,k)$ as a Banach $k$-algebra. In the case where $X$ is compact, the underlying $k$-vector space of $\t{C}_{\t{bd}}(X,k)$ coincides with $\t{C}_0(X,k)$, and hence the identification above involves no ambiguity.

\begin{exm}
\label{Tate algebra}
Let $r \in (0,\infty)$. The $k$-algebra
\begin{eqnarray*}
  k \Ens{r^{-1}T} \coloneqq \Set{F = \sum_{n = 0}^{\infty} F_n T^n \in k[[T]]}{\lim_{n \to \infty} \vv{F_n} r^n = 0}
\end{eqnarray*}
is a commutative uniform Banach $k$-algebra with respect to the multiplicative norm
\begin{eqnarray*}
  \nn{\cdot} \colon k \Ens{r^{-1}T} & \to & [0,\infty) \\
  \sum_{n = 0}^{\infty} F_n T^n & \mapsto & \nn{F} \coloneqq \sup_{n \in \N} \vv{F_n} r^n < \infty
\end{eqnarray*}
called {\it the Gauss norm of radius $r$}. The evaluation map
\begin{eqnarray*}
  k \Ens{r^{-1}T} & \hookrightarrow & \t{C}_{\t{bd}}(k(r),k) \\
  \sum_{n = 0}^{\infty} F_n T^n & \mapsto & \left( \lambda \mapsto \sum_{n = 0}^{\infty} F_n \lambda^n \right)
\end{eqnarray*}
is a submetric injective $k$-linear homomorphism. We identify the underlying $k$-algebra of $k \ens{r^{-1}T}$ with its image $k \ens{r^{-1} z} \subset \t{C}_{\t{bd}}(k(r),k)$. When $r = 1$, then we abbreviate $k \ens{r^{-1}T}$ by $k \ens{T}$. The natural embedding $k[T] \hookrightarrow k \ens{T}$ induces a $\overline{k}$-algebra isomorphism $\overline{k \ens{T}} \to \overline{k}[T]$.
\end{exm}

\begin{exm}
Let $(r_1,r_2) \in (0,\infty) \times (0,\infty)$ with $r_1 \leq r_2$. The $k$-algebra
\begin{eqnarray*}
  k \Ens{r_2^{-1} T, r_1 T^{-1}} \coloneqq \Set{F = \sum_{n = - \infty}^{\infty} F_n T^n \in \prod_{n \in \Z} k T^n}{\lim_{n \to \infty} \max \Ens{\vv{F_n} r_1^n,  \vv{F_{-n}} r_2^{-n}} = 0}
\end{eqnarray*}
is a commutative uniform Banach $k$-algebra with respect to the power-multiplicative norm
\begin{eqnarray*}
  \nn{\cdot} \colon k \Ens{r_2^{-1} T, r_1 T^{-1}} & \to & [0,\infty) \\
  \sum_{n = - \infty}^{\infty} F_n T^n & \mapsto & \nn{F} \coloneqq \sup_{n \in \N} \max \Ens{\vv{F_n} r_1^n,  \vv{F_{-n}} r_2^{-n}} < \infty
\end{eqnarray*}
called {\it the Gauss norm of radius $(r_1,r_2)$}. The evaluation map
\begin{eqnarray*}
  k \Ens{r_2^{-1}T, r_1 T^{-1}} & \hookrightarrow & \t{C}_{\t{bd}}(k(r_2) \backslash k(r_1-), k) \\
  \sum_{n = - \infty}^{\infty} F_n T^n & \mapsto & \left( \lambda \mapsto \sum_{n = - \infty}^{\infty} F_n \lambda^n \right)
\end{eqnarray*}
is a submetric injective $k$-linear homomorphism. We identify the underlying $k$-algebra of $k \ens{r_2^{-1} T, r_1 T^{-1}}$ with its image $k \ens{r_2^{-1} z, r_1 z^{-1}} \subset \t{C}_{\t{bd}}(k(r_2) \backslash k(r_1-),k)$. When $r_1 = r_2 = 1$, then we abbreviate $k \ens{r_2^{-1}T, r_1T^{-1}}$ by $k \ens{T,T^{-1}}$. The natural embedding $k[T,T^{-1}] \hookrightarrow k \ens{T,T^{-1}}$ induces a $\overline{k}$-algebra isomorphism $\overline{k \ens{T,T^{-1}}} \to \overline{k}[T,T^{-1}]$.
\end{exm}

\begin{dfn}
Let $\A$ be a Banach $k$-algebra. For each $A \in \A$, we denote by $k[A] \subset \A$ the image of the $k$-algebra homomorphism $k[T] \to \A \colon T \mapsto A$, by $\Leb(A) \subset \A$ the closure of $k[A]$, and by $\Leb_{\A}(A) \subset \A$ the closure of the localisation of $k[A]$ by the multiplicative subset $k[A] \cap \A^{\times}$.
\end{dfn}

The $k$-subalgebra $\Leb(A) \subset \A$ is a commutative Banach $k$-algebra, whose isomorphism class is independent of the choice of $\A$. Namely, even if one replaces $\A$ by a closed $k$-subalgebra $\B \subset \A$ containing $A$, $\Leb(A) \subset \B$ does not change as a $k$-subalgebra of $\A$. On the other hand, $\Leb_{\A}(A)$ is a commutative Banach $k$-algebra whose isomorphism class depends on the choice of the Banach $k$-algebra $\A$ containing $A$.

\begin{exm}
$\t{ }$
\begin{itemize}
\item[(i)] For any $r \in (0,\infty)$, we have
\begin{eqnarray*}
  \Leb(T) = \Leb_{k \Ens{r^{-1} T}}(T) = k \Ens{r^{-1}T}
\end{eqnarray*}
\item[(ii)] For any $(r_1,r_2) \in (0,\infty)$ with $r_1 \leq r_2$, we have
\begin{eqnarray*}
  \Leb(T) = k \Ens{r_2^{-1} T} \subsetneq \Leb_{k \Ens{r_2^{-1} T, r_1 T^{-1}}}(T) = k \Ens{r_2^{-1} T, r_1 T^{-1}}.
\end{eqnarray*}
\end{itemize}
\end{exm}

\subsection{Rigidity of Continuous Functions}
\label{Rigidity of Continuous Functions}

In this subsection, we deal with several properties of bounded continuous functions. We will use the results in \S \ref{Rigid Continuous Functional Calculus} in order to formulate and study functional calculi.

\begin{prp}
\label{reduction of continuous functions}
Let $X$ be a topological space. If $X$ is compact, or if the valuation of $k$ is discrete or trivial, then, there is a natural $\overline{k}$-algebra isomorphism $\overline{\t{C}_{\t{bd}}(X,k)} \to \t{C}(X,\overline{k})$.
\end{prp}

\begin{proof}
We identify $\t{C}_{\t{bd}}(X,k)(1)$ and $\t{C}_{\t{bd}}(X,O_k)$ through the inclusion $O_k \hookrightarrow k$. The kernel of the $O_k$-algebra homomorphism $\t{C}_{\t{bd}}(X,k)(1) \to \t{C}(X,\overline{k})$ induced by the reduction $O_k \twoheadrightarrow \overline{k}$ coincides with $\t{C}_{\t{bd}}(X,k)(1-)$ by the maximum modulus principle. Therefore it induces an injective $\overline{k}$-algebra homomorphism $\iota \colon \overline{\t{C}_{\t{bd}}(X,k)} \to \t{C}(X,\overline{k})$. Let $\lbrack \cdot \rbrack \colon \overline{k} \hookrightarrow k$ be a Teichm\"uller embedding, which is continuous because of the discreteness of the topology of $\overline{k}$. The composite with $\lbrack \cdot \rbrack$ gives a section of $\iota$. It implies that $\iota$ is surjective.
\end{proof}

In the following, let $\sigma \subset k$ denote a non-empty bounded closed subset, $z_{\sigma}^{(k)} \colon \sigma \hookrightarrow k$ the coordinate function given as the inclusion. The boundedness of $\sigma$ ensures the boundedness of $z_{\sigma}$. We abbreviate $z_{\sigma}^{(k)}$ to $z_{\sigma}$ as long as there is no ambiguity of $k$, and $z_{\sigma}$ to $z$ as long as there is no ambiguity of $\sigma$.

\begin{dfn}
We set $\t{C}_{\t{rig}}(\sigma,k) \coloneqq \Leb_{\t{C}_{\t{bd}}(\sigma,k)}(z) \subset \t{C}_{\t{bd}}(\sigma,k)$, and call an element of $\t{C}_{\t{rig}}(\sigma,k)$ a {\it rigid continuous function on $\sigma$}.
\end{dfn}

A reader might replace the term ``rigid continuous'' by ``rigid analytic'' or ``Krasner analytic''. We aboid from use of those terms just because they are quite ambiguous.

\begin{exm}
\label{rigidity}
$\t{ }$
\begin{itemize}
\item[(i)] In the case $k = \Q_p$ and $\sigma = \Z_p \subset \Q_p$, we have
\begin{eqnarray*}
  \Leb(z) = \t{C}_{\t{rig}}(\Z_p,\Q_p) = \t{C}_{\t{bd}}(\Z_p,\Q_p) = \t{C}(\Z_p,\Q_p).
\end{eqnarray*}
\item[(ii)] In the case $k = \Q_p$ and $\sigma = \Z_p \backslash p \Z_p \subset \Q_p$, we have
\begin{eqnarray*}
  \Leb(z) = \t{C}_{\t{rig}}(\Z_p \backslash p \Z_p, \Q_p) = \t{C}_{\t{bd}}(\Z_p \backslash p \Z_p, \Q_p) = \t{C}(\Z_p \backslash p \Z_p, \Q_p).
\end{eqnarray*}
\item[(iii)] In the case $k = \C_p$ and $\sigma = O_{\C_p} \subset \C_p$, we have
\begin{eqnarray*}
  \Leb(z) = \t{C}_{\t{rig}}(O_{\C_p},\C_p) = \C_p \Ens{z} \subsetneq \t{C}_{\t{bd}}(O_{\C_p},\C_p).
\end{eqnarray*}
\item[(iv)] In the case $k = \C_p$ and $\sigma = O_{\C_p} \backslash p m_{\C_p} \subset \C_p$, we have
\begin{eqnarray*}
  \Leb(z) \subsetneq \t{C}_{\t{rig}}(O_{\C_p} \backslash p m_{\C_p}, \C_p) = \C_p \Ens{z, \vv{p} z^{-1}} \subsetneq \t{C}_{\t{bd}}(O_{\C_p} \backslash p m_{\C_p}, \C_p).
\end{eqnarray*}
\end{itemize}
\end{exm}

\begin{prp}
\label{continuous functions on a compact subset}
If $\sigma$ is compact, then $k[z]$ is dense in $\t{C}_{\t{rig}}(\sigma,k)$ and the equality $\Leb(z) = \t{C}_{\t{rig}}(\sigma,k) = \t{C}(\sigma,k)$ holds.
\end{prp}

\begin{proof}
The assertion is verified in \cite{Mur78} 1.7 and \cite{Ber90} 9.2.6. We give a brief explanation of the proof. There is a well-known $p$-adic analogue of Stone--Weierstrass theorem, which is originally verified by Kaplansky in \cite{Kap50} as a generalisation of Dieudonn\'e's theorem in \cite{Die44}. It states that $k[z]$ is dense in $\t{C}(\sigma,k)$. Therefore the equality in the assertion holds because we have $k[z] \subset \Leb(z) \subset \t{C}_{\t{rig}}(\sigma,k) \subset \t{C}(\sigma,k)$.
\end{proof}

\begin{lmm}
\label{invertible polynomial function}
For any $P \in k[T]$ with no zeros on $\sigma$, the polynomial function
\begin{eqnarray*}
  P(z) \colon \sigma & \to & k \\
  \lambda & \mapsto & P(\lambda)
\end{eqnarray*}
is invertible in $\t{C}_{\t{rig}}(\sigma,k)$.
\end{lmm}

\begin{proof}
By the definition of $\t{C}_{\t{rig}}(\sigma,k)$, it suffices to show that $P(z)$ is invertible in $\t{C}_{\t{bd}}(\sigma,k)$. Take an algebraic closure $k^{\t{alg}}$ of $k$, and equip $k^{\t{alg}}$ with the ultrametric given by a unique extension of the valuation of $k$. We denote by $S \subset k^{\t{alg}} \backslash \sigma$ the finite set of the zeros of $P$. Since $S$ is compact and $\sigma$ is closed, the distance $d$ between $S$ and $\sigma$ is positive. Put $P = c \prod_{s \in S} (T - s)^{n_s}$ by a $(c,(n_s)_{s \in S}) \in k^{\times} \times (\N \backslash \ens{0})^S$. Set $n \coloneqq \sum_{s \in S} n_s$. We have $\v{P(\lambda)} = \v{c} \prod_{s \in S} \v{\lambda - s}^{n_s} \geq \v{c} d^n > 0$ for any $\lambda \in \sigma$, and hence the real-valued function $\v{P(z)} \colon \sigma \to [0,\infty) \colon \lambda \mapsto \v{P(\lambda)}$ admits a positive lower bound $\v{c} d^n$. It implies that $P(z)$ is invertible in $\t{C}_{\t{bd}}(\sigma,k)$.
\end{proof}

\begin{thm}
\label{invertible rigid continuous function}
The equalities $\t{C}_{\t{rig}}(\sigma,k)^{\times} = \t{C}_{\t{rig}}(\sigma,k) \cap \t{C}_{\t{bd}}(\sigma,k)^{\times}$ and $\sigma_{\t{C}_{\t{rig}}(\sigma,k)}(z_{\sigma}) = \sigma$ hold.
\end{thm}

\begin{proof}
We have $\t{C}_{\t{rig}}(\sigma,k)^{\times} \subset \t{C}_{\t{rig}}(\sigma,k) \cap \t{C}_{\t{bd}}(\sigma,k)^{\times}$. Take an $f \in \t{C}_{\t{rig}}(\sigma,k) \cap \t{C}_{\t{bd}}(\sigma,k)^{\times}$. By the definition of $\t{C}_{\t{rig}}(\sigma,k)$, for any $\epsilon \in (0,1)$, there is a continuous function $g \colon \sigma \to k$ contained in the localisation of $k[z] \subset \t{C}_{\t{bd}}(\sigma,k)$ with $\n{f-g} < \epsilon \min \ens{\n{f^{-1}}^{-1}, \n{f^{-1}}^{-2}}$. For any $\lambda \in \sigma$, we have $\n{f(\lambda)} = \n{f^{-1}(\lambda)}^{-1} \geq \n{f^{-1}}^{-1} > \n{f-g} \geq \v{f(\lambda) - g(\lambda)}$, and hence $\v{g(\lambda)} = \v{f(\lambda)} > 0$. Therefore we obtain $\v{g(\lambda)}^{-1} = \v{f(\lambda)}^{-1} = \v{f^{-1}(\lambda)} \leq \n{f^{-1}}$. It implies that $g$ is invertible in $\t{C}_{\t{bd}}(\sigma,k)$. Take a $(P,Q) \in k[T] \times k[T]$ with $Q(z) \in \t{C}_{\t{bd}}(\sigma,k)^{\times}$ and $g = Q(z)^{-1}P(z)$. Since $Q(z)$ and $g$ belong to $\t{C}_{\t{bd}}(\sigma,k)^{\times}$, $Q(z)$ and $g$ have no zeros on $\sigma$, and hence so does $P(z) = gQ(z)$. Therefore $P(z)$ is invertible in $\t{C}_{\t{bd}}(\sigma,k)$ by Lemma \ref{invertible polynomial function}, and $g^{-1} = P(z)^{-1}Q(z)$ is contained in the localisation of $k[z]$. We have $\n{f^{-1} - g^{-1}} \leq \n{f^{-1}} \ \n{g^{-1}} \ \n{g - f} = \n{f^{-1}}^2 \n{f - g} < \epsilon$, and hence $f \in \t{C}_{\t{rig}}(\sigma,k)$. It implies $\t{C}_{\t{rig}}(\sigma,k)^{\times} = \t{C}_{\t{rig}}(\sigma,k) \cap \t{C}_{\t{bd}}(\sigma,k)^{\times}$. We have $\sigma_{\t{C}_{\t{rig}}(\sigma,k)}(z) \subset \sigma$ because $\sigma$ coincides with the image of $z$, and also $\sigma \subset \sigma_{\t{C}_{\t{rig}}(\sigma,k)}(z)$ by the definition of $\t{C}_{\t{rig}}(\sigma,k)$. Therefore $\sigma_{\t{C}_{\t{rig}}(\sigma,k)}(z)$ coincides with $\sigma$.
\end{proof}

\begin{crl}
\label{invertible rigid continuous function 2}
If $\sigma$ is contained in $O_k$, then the equality $\sigma_{\t{C}_{\t{rig}}(\sigma,k)(1)}(z_{\sigma}) = \set{\lambda + \varpi}{(\lambda,\varpi) \in \sigma \times m_k}$ holds.
\end{crl}

\begin{proof}
Let $(\lambda,\varpi) \in \sigma \times m_k$. Then $z_{\sigma} - (\lambda + \varpi)$ sends $\lambda \in \sigma$ to $\varpi \in m_k$, and hence is not invertible in $\t{C}_{\t{bd}}(\sigma,k)(1)$. Therefore $z_{\sigma} - (\lambda + \varpi)$ is not invertible in $\t{C}_{\t{rig}}(\sigma,k)(1)$, and $\lambda$ lies in $\sigma_{\t{C}_{\t{rig}}(\sigma,k)(1)}(z_{\sigma})$. It implies $\set{\lambda + \varpi}{(\lambda,\varpi) \in \sigma \times m_k} \subset \sigma_{\t{C}_{\t{rig}}(\sigma,k)(1)}(z_{\sigma})$.

\vspace{0.2in}
Let $\lambda \in k \backslash \set{\lambda + \varpi}{(\lambda,\varpi) \in \sigma \times m_k}$. Then $z_{\sigma} - \lambda$ factors through the inclusion $O_k^{\times} \hookrightarrow k$, and hence is invertible in $\t{C}_{\t{rig}}(\sigma,k)(1)$ by Theorem \ref{invertible rigid continuous function}. Therefore $\lambda$ lies in $k \backslash \sigma_{\t{C}_{\t{rig}}(\sigma,k)(1)}(z_{\sigma})$. It implies $\sigma_{\t{C}_{\t{rig}}(\sigma,k)(1)}(z_{\sigma}) = \set{\lambda + \varpi}{(\lambda,\varpi) \in \sigma \times m_k}$.
\end{proof}

We show certain functoriality of the correspondence $\sigma \rightsquigarrow \t{C}_{\t{rig}}(\sigma,k)$.

\begin{dfn}
Let $X$ be a topological space, and $X_0 \subset X$ a clopen subspace. We denote by $\delta^{(k)}_{X_0,X} \colon X \to k$ the characteristic function of $X_0$. We abbreviate $\delta^{(k)}_{X_0,X}$ to $\delta_{X_0,X}$ as long as there is no ambiguity of $k$. When $X$ is a discrete topological space $I$ and $X_0$ is a singleton $\ens{i}$ for an $i \in I$, then we abbreviate $\delta^{(k)}_{\ens{i},I}$ to $\delta^{(k)}_i$ as long as there is no ambiguity of $I$.
\end{dfn}

\begin{rmk}
Let $I$ be a set. The collection $\set{\delta_i}{i \in I}$ forms an orthonormal Schauder basis of $\t{C}_0(I,k)$, but does not form an orthonormal Schauder basis of $\t{C}_{\t{bd}}(I,k)$ unless $I$ is a finite set.
\end{rmk}

\begin{dfn}
\label{measurable}
A clopen subset $\sigma_0 \subset \sigma$ is said to be {\it measurable} if $\delta_{\sigma_0,\sigma} \colon \sigma \to k$ lies in $\t{C}_{\t{rig}}(\sigma,k)$.
\end{dfn}

\begin{exm}
If $\sigma \subset k$ is compact, then every clopen subset of $\sigma$ is measurable by Proposition \ref{continuous functions on a compact subset}. 
\end{exm}

\begin{thm}
\label{functoriality}
Let $\sigma_0 \subset \sigma$ be a closed subset. Then the image of $\t{C}_{\t{rig}}(\sigma,k)$ by the restriction map $\t{C}_{\t{bd}}(\sigma,k) \twoheadrightarrow \t{C}_{\t{bd}}(\sigma_0,k)$ is contained in $\t{C}_{\t{rig}}(\sigma_0,k)$. When $\sigma_0$ is a measurable subset of $\sigma$, then the image of $\t{C}_{\t{rig}}(\sigma,k)$ by the restriction map $\t{C}_{\t{bd}}(\sigma,k) \twoheadrightarrow \t{C}_{\t{bd}}(\sigma_0,k)$ coincides with $\t{C}_{\t{rig}}(\sigma_0,k)$, and the image of $\t{C}_{\t{rig}}(\sigma_0,k)$ by the zero-extension $\t{C}_{\t{bd}}(\sigma_0,k) \hookrightarrow \t{C}_{\t{bd}}(\sigma,k)$ coincides with $\delta_{\sigma_0,\sigma} \t{C}_{\t{rig}}(\sigma,k)$.
\end{thm}

\begin{proof}
Let $\Pi \colon \t{C}_{\t{bd}}(\sigma,k) \twoheadrightarrow \t{C}_{\t{bd}}(\sigma_0,k)$ denote the restriction map. The first assertion follows from the continuity of $\Pi$. Suppose that $\sigma_0$ is a measurable subset of $\sigma$. The image of $\t{C}_{\t{bd}}(\sigma_0,k)$ by the zero-extension $\t{C}_{\t{bd}}(\sigma_0,k) \hookrightarrow \t{C}_{\t{bd}}(\sigma,k)$ coincides with $\delta_{\sigma_0,\sigma} \t{C}_{\t{bd}}(\sigma_0,k)$. We denote by $\iota_! \colon \t{C}_{\t{bd}}(\sigma_0,k) \to \delta_{\sigma_0,\sigma} \t{C}_{\t{bd}}(\sigma_0,k)$ the multiplicative $k$-linear isomorphism induced by the zero-extension, which is an isometry by the definition of the supremum norm. The restriction of $\Pi$ on $\delta_{\sigma_0,\sigma} \t{C}_{\t{bd}}(\sigma_0,k)$ is a multiplicative $k$-linear isomorphism, and is the inverse of $\iota_!$. For any $f \in \t{C}_{\t{bd}}(\sigma_0,k)^{\times}$, $\iota_!(f) + (1-\delta_{\sigma_0,\sigma})$ admits the inverse $\iota_!(f^{-1}) + (1-\delta_{\sigma_0,\sigma})$ in $\t{C}_{\t{bd}}(\sigma_0,k)$, and the equality $\Pi(\iota_!(f) + (1-\delta_{\sigma_0,\sigma})) = f$. Therefore the image of $\t{C}_{\t{bd}}(\sigma,k)^{\times}$ by $\Pi$ coincides with $\t{C}_{\t{bd}}(\sigma_0,k)^{\times}$.

\vspace{0.2in}
Let $f \in \t{C}_{\t{rig}}(\sigma_0,k)$. By the definition of $\t{C}_{\t{rig}}(\sigma_0,k)$, $f$ is the limit of a sequence on the localisation of $k[z_{\sigma_0}]$ by $k[z_{\sigma_0}] \cap \t{C}_{\t{bd}}(\sigma_0,k)^{\times}$. By the equalities $z_{\sigma_0} = \Pi(z_{\sigma})$ and $\iota_!(z_{\sigma_0}) = \delta_{\sigma_0,\sigma} z_{\sigma}$, $\iota_!(f)$ is the limit of a sequence on the localisation of $k[\delta_{\sigma_0,\sigma} z_{\sigma}]$ by $k[\delta_{\sigma_0,\sigma} z_{\sigma}] \cap \t{C}_{\t{bd}}(\sigma,k)^{\times}$. Since $\sigma_0$ is a measurable subset of $\sigma$, $\iota_!(f)$ lies in $\t{C}_{\t{rig}}(\sigma,k)$. It implies $\iota_!(\t{C}_{\t{rig}}(\sigma_0,k)) \subset \t{C}_{\t{rig}}(\sigma,k)$ and hence $\t{C}_{\t{rig}}(\sigma_0,k) \subset (\Pi \circ \iota_!)(\t{C}_{\t{rig}}(\sigma_0,k)) \subset \Pi(\t{C}_{\t{rig}}(\sigma,k))$. Thus we obtain $\t{C}_{\t{rig}}(\sigma_0,k) = \Pi(\t{C}_{\t{rig}}(\sigma,k))$.
\end{proof}

\subsection{Operator Algebra}
\label{Operator Algebra}

In this subsection, let $V$ denote a Banach $k$-vector spaces. A {\it bounded operator on $V$} is a bounded $k$-linear endomorphism $f \colon V \to V$. We denote by $\B_k(V)$ the $k$-algebra of bounded operators on $V$, and endow it with the norm $\n{\cdot} \colon \B_k(V) \to [0,\infty)$ called {\it the operator norm} given by setting $\sup_{v \in V(1)} \n{Av} < \infty$ for each $A \in \B_k(V)$. The completeness of $V$ ensures that $\B_k(V)$ is a Banach $k$-algebra.

\begin{exm}
Let $n \in \N$. Then the canonical basis of $k^n$ induces a bijective $k$-linear homomorphism $\B_k(k^n) \to \t{M}_n(k)$, and hence we identify $\B_k(k^n)$ and $\t{M}_n(k)$.
\end{exm}

\begin{exm}
If the valuation of $k$ is trivial and $V$ is a discrete $k$-vector space, then $\B_k(V)$ is the discrete $k$-algebra $\t{End}_k(V)$ of $k$-linear endomorphisms.
\end{exm}

\begin{dfn}
Let $I$ be a set. For an $A \in \B_k(\t{C}_0(I,k))$, {\it the matrix representation of $A$ with respect to $I$} is the function $M_{I,k}(A) \colon I \times I \to k \colon (i,j) \mapsto M_{I,k}(A)_{i,j}$ given by setting $M_{I,k}(A)_{i,j} \colon (A \delta_i)(j)$ for each $(i,j) \in I \times I$.
\end{dfn}

\begin{prp}
\label{matrix representation}
Let $I$ be a set. Then $M_{I,k}(A)$ is a bounded function for any $A \in \B_k(\t{C}_0(I,k))$, and the well-defined map $M_{I,k} \colon \B_k(\t{C}_0(I,k)) \to \t{C}_{\t{bd}}(I \times I,k) \colon A \mapsto M_{I,k}(A)$ is an isometric $k$-linear isomorphism onto the closed image.
\end{prp}

\begin{proof}
Let $A \in \B_k(\t{C}_0(I,k))$. For any $(i,j) \in I \times I$, we have $\v{M_{I,k}(A)_{i,j}} \leq \n{A \delta_i} \leq \n{A}$. Therefore $M_{I,k}(A)$ is a bounded function with $\n{M_{I,k}(A)} \leq \n{A}$. Let $f \in \B_k(\t{C}_0(I,k)(1))$. By the definition of $\t{C}_0(I,k)$, $f(i) = 0$ all but countably many $i \in I$, and the essentially countable infinite sum $\sum_{i \in I} f(i) \delta_i$ converges to $f$. More precisely, $(\sum_{i \in F} f(i) \delta_i)_{F \in \Fil(I)}$ is a net on $\t{C}_0(I,k)$ converging to $f$, where $\Fil(I)$ is the set of finite subsets of $I$ directed by inclusions. By the continuity of $A$, $(\sum_{i \in F} f(i) A \delta_i)_{F \in \Fil(I)}$ is a net on $\t{C}_0(I,k)$ converging to $Af$. Therefore we obtain $\n{Af} \leq \sup_{F \in \Fil(I)} \max_{i \in F} \v{f(i)} \n{A \delta_i} \leq \n{f} \n{M_{I,k}(A)} \leq \n{M_{I,k}(A)}$. It implies $\n{A} \leq \n{M_{I,k}(A)}$. Thus we obtain a well-defined isometric $k$-linear homomorphism $M_{I,k} \colon \B_k(\t{C}_0(I,k)) \to \t{C}_{\t{bd}}(I \times I,k) \colon A \mapsto M_{I,k}(A)$, whose image is closed because of the completeness of $\B_k(\t{C}_0(I,k))$.
\end{proof}

For a set $I$, we abbreviate $M_{I,k}$ to $M_I$ as long as there is no ambiguity of $k$. The notion of the matrix representation of a bounded operator is a generalisation of the matrix representation of a $k$-linear endomorphism on a finite dimensional $k$-vector space.

\begin{exm}
Let $n \in \N$. For any $A \in \B_k(k^{\N \cap [1,n]})$, $M_{\N \cap [1,n]}(A)$ coincides with the image of $A$ by the identification $\B_k(k^{\N \cap [1,n]}) \to \B_k(k^n) = \t{M}_n(k)$ induced by the canonical basis $(\delta_i)_{i \in \N \cap [1,n]}$ of $k^{\N \cap [1,n]}$.
\end{exm}

\begin{exm}
\label{reductive matrix representation}
Suppose that the valuation of $k$ is trivial. Let $I$ be a set. For any $A \in \t{End}_k(k^{\oplus I})$, we have $(A \delta_i(j))_{j \in I} \in k^{\oplus I}$ for any $i \in I$, and the image of $M_I \colon \t{End}_k(k^{\oplus I}) \to k^{I \times I}$ coincides with the image of the natural embedding $(k^{\oplus I})^I \hookrightarrow (k^I)^I \cong k^{I \times I}$. Therefore we regard $M_I$ as an isometric $k$-linear isomorphism $\t{End}_k(k^{\oplus I}) \to (k^{\oplus I})^I$ with respect to the trivial norm on $(k^{\oplus})^I$.
\end{exm}

\begin{dfn}
Suppose that $V$ is a strictly Cartesian $k$-vector space. For an $A \in \B_k(V)(1)$, {\it the operator reduction of $A$} is the image $\Pi_V(A)$ of $A$ by the continuous surjective $O_k$-algebra homomorphism $\Pi_V \colon \B_k(V)(1) \twoheadrightarrow \t{End}_{\overline{k}}(\overline{V})$ given by the functoriality of the reduction $V \rightsquigarrow \overline{V}$ with respect to submetric $k$-linear homomorphism. A bounded operator $A$ on $V$ is said to be {\it reductively scalar} if $A$ is submetric and $\Pi_V(A)$ lies in the image of $\overline{k}$, and is said to be {\it reductively transcendental} if $A$ is submetric and $\Pi_V(A)$ is not integral over $\overline{k}$.
\end{dfn}

\begin{exm}
Let $n \in \N$. We identify $\t{M}_n(k)$ with $\B_k(k^n)$ in a natural way, where $k^n$ is endowed with the norm associated to the canonical basis. Let $M \in \t{M}_n(k)(1) = \t{M}_n(O_k)$. Then $\Pi_{k^n}(M)$ is the element of $\t{End}_{\overline{k}}(\overline{k}^n) = \t{M}_n(\overline{k})$ whose entries are the reductions of the entries of $M$. Therefore $M$ is reductively scalar if and only if every entry other than diagonal entries lies in $m_k$, and $M$ is never reductively transcendental. We note that a replacement of the norm of $k^n$ by an equivalent one for which $k^n$ is strictly Cartesian possibly causes a change of the integral model $\t{M}_n(k^n)(1) \subset \t{M}_n(k^n)$, and hence of the matrix reduction.
\end{exm}

The following two propositions ensure that the operator reduction of an operator is a generalisation of a matrix, which is given as the reduction of each entry, and is not identified with the image by the reduction of the operator algebra.

\begin{prp}
\label{operator reduction & matrix representation}
Let $I$ be a set, and $A \in \B_k(\t{C}_0(I,k))(1)$. Then the equality
\begin{eqnarray*}
  M_{I,\overline{k}}(\Pi_{\t{C}_0(I,k)}(A))_{i,j} = M_{I,k}(A)_{i,j} + m_k
\end{eqnarray*}
holds for any $(i,j) \in I \times I$.
\end{prp}

\begin{proof}
For any $(i,j) \in I \times I$, we have
\begin{eqnarray*}
  & & M_{I,\overline{k}} \left( \Pi_{\t{C}_0(I,k)}(A) \right)_{i,j} = \left( \Pi_{\t{C}_0(I,k)}(A) \delta^{(\overline{k})}_i \right)(j) = \left( A \delta^{(k)}_i + \t{C}_0(I,k)(1-) \right)(j) \\
  & = & \left( A \delta^{(k)}_i \right)(j) + m_k = M_{I,k}(A)_{i,j} + m_k
\end{eqnarray*}
through the identification $\overline{\t{C}_0(I,k)} \cong \overline{k}{}^{\oplus I}$ in Proposition \ref{reduction & direct product}. Thus the assertion holds.
\end{proof}

\begin{prp}
\label{reduction & operator reduction}
Suppose that $V$ is a strictly Cartesian $k$-vector space. Then the surjective $\overline{k}$-linear homomorphism $\overline{\Pi}_V \colon \overline{\B_k(V)} \twoheadrightarrow \t{End}_{\overline{k}}(\overline{V})$ induced by $\Pi_V$ is an isomorphism if and only if $V$ is of finite dimension or the valuation of $k$ is discrete or trivial.
\end{prp}

\begin{proof}
By the functoriality of the reduction and the naturalness of $\Pi_V$, we may assume $V = \t{C}_0(I,k)$ for a set $I$ without loss of generality. Suppose that $I$ is a finite set or the valuation of $k$ is discrete or trivial. In order to show the injectivity of $\overline{\Pi}_{\t{C}_0(I,k)}$, it suffices to verify that every $A \in \B_k(\t{C}_0(I,k))(1)$ with $\Pi_{\t{C}_0(I,k)}(A) = 0$ lies in $\B_k(\t{C}_0(I,k))(1-)$. Let $A \in \B_k(\t{C}_0(I,k))(1)$ with $\Pi_{\t{C}_0(I,k)}(A) = 0$. By Proposition \ref{operator reduction & matrix representation}, we have $M_{I,k}(A)_{i,j} + m_k = M_{I,\overline{k}}(\Pi_{\t{C}_0(I,k)}(A))_{i,j} = 0$ for any $(i,j) \in I \times I$, and hence $M_{I,k}(A)$ gives a function $I \times I \to m_k$. If $I$ is a finite set, then so is the image of $M_{I,k}(A)$, and hence there is an $r \in (0,1)$ with $\v{M_{I,k}(A)_{i,j}} \leq r$ for any $(i,j) \in I \times I$. If the valuation of $k$ is discrete or trivial, then there is an $r \in (0,1)$ with $m_k \subset k(r)$, and we have $\v{M_{I,k}(A)_{i,j}} \leq r$ for any $(i,j) \in I \times I$. Therefore Proposition \ref{matrix representation} ensures $\n{A} \leq r < 1$. We obtain $A \in \B_k(\t{C}_0(I,k))(1-)$. Thus $\Pi_{\t{C}_0(I,k)}$ is injective.

\vspace{0.2in}
Suppose that $I$ is an infinite set and the valuation of $k$ is neither discrete nor trivial. Take an injective map $i_{\bullet} \colon \N \hookrightarrow I \colon n \mapsto i_n$ and a sequence $(\varpi_n)_{n \in \N} \in k^{\N}$ such that $(\v{\varpi_n})_{n \in \N}$ is a strictly increasing sequence converging to $1$. Let $A$ denote a bounded operator on $\t{C}_0(I,k)$ given by setting $(Af) = (\sum_{n = 0}^{\infty} \varpi_n f(i_n)) \delta^{(k)}_{i_0}$ for any $f \in \t{C}_0(I,k)$. Then we have $\n{A \delta^{(k)}_{i_n}} = \n{\varpi_n \delta_{i_0}} = \v{\varpi_n}$ for any $n \in \N$ and $\n{A \delta^{(k)}_i} = 0$ for any $i \in I \backslash \set{i_n}{n \in \N}$. It implies $\n{A} = 1$. On the other hand, we have $\Pi_{\t{C}_0(I,k)}(A) \delta^{(\overline{k})}_{i_n} = \varpi_n f(i_n) + m_k = 0$ for any $n \in \N$ and $\Pi_{\t{C}_0(I,k)}(A) \delta^{(\overline{k})}_i = 0$ for any $i \in I \backslash \set{i_n}{n \in \N}$. It implies $\Pi_{\t{C}_0(I,k)}(A) = 0$. Thus $\overline{\Pi}_{\t{C}_0(I,k)}$ is not injective.
\end{proof}

Let $V$ be a Banach $k$-vector space. A {\it unitary operator on $V$} is an isometric $k$-linear automorphism on $V$. In particular, every unitary operator on $V$ is a bounded operator $V$, and a $U \in \B_k(V)$ is a unitary operator if and only if $U$ lies in $\B_k(V)(1)^{\times}$.

\begin{prp}
\label{criterion for unitarity}
Let $n \in \N$ and $U \in \t{M}_n(k)$. Then the following are equivalent:
\begin{itemize}
\item[(i)] The operator $U$ is unitary.
\item[(ii)] The collection of columns of $U$ forms an orthonormal basis of $k^n$.
\item[(iii)] The operator $U$ lies in $\t{M}_n(O_k)$ and $\Pi_{k^n}(U)$ lies in $\t{GL}_n(\overline{k})$.
\item[(iv)] The equalities $\n{U} = 1$ and $\v{\det U} = 1$ hold.
\item[(v)] The operator $U$ lies in $\t{GL}_n(O_k)$.
\item[(vi)] The operator $U$ lies in $\t{GL}_n(k)$ and $U^{-1}$ is unitary.
\item[(vii)] The operator ${}^{\t{t}} U$ is unitary.
\end{itemize}
\end{prp}

\begin{proof}
Assume (i). Since $U$ is isometric, we have $\n{\sum_{i = 1}^{n} c_i U \delta_i} = \n{U(\sum_{i = 1}^{n} c_i \delta_i)} = \n{\sum_{i = 1}^{n} c_i \delta_i} = \max_{1 \leq i \leq n} \v{c_i}$ for any $(c_i)_{i = 1}^{n} \in k^n$. It implies that $(U \delta_i)_{i = 1}^{n}$ forms an orthonormal basis of $k^n$. Therefore $U$ satisfies (ii). Assume (ii). Since every element of an orthonormal basis of $k^n$ lies in $O_k^n$, we have $U \in \t{M}_n(O_k)$. Moreover, the reduction of an orthonormal basis of $k^n$ forms a basis in $\overline{k}{}^n$ by \cite{BGR84} 2.5.1/4. Therefore $U$ satisfies (iii). Assume (iii). We have $(\det U) + m_k = \det \Pi_{k^n}(U) \neq 0$ by Proposition \ref{operator reduction & matrix representation}. It implies $\det U \notin m_k$ and hence $U \notin m_k \t{M}_n(O_k)$. Therefore $U$ satisfies (iv). Assume (iv). We have $U \in \t{M}_n(O_k)$ and $\det U \in O_k^{\times}$. It implies that $U$ satisfies (v). Assume (v). We have $U \in \t{GL}_n(k)$ and $U^{-1} \in \t{GL}_n(O_k)$. It implies that $U$ and $U^{-1}$ preserve the closed unit ball $O_k^n \subset k^n$, and hence $U$ satisfies (vi). Since (i) for $U$ implies (vi) for $U$, (vi) for $U$ implies (vi) for $U^{-1}$ and hence (i) for $U$. Since (v) for $U$ is equivalent to (v) for ${}^{\t{t}} U$, (vii) for $U$ is equivalent to (i) for $U$.
\end{proof}

The notion of the unitarity of a matrix deeply depends on the choice of the norm of $k^n$. Let $U \in \t{GL}_n(k)$ be a diagonalisable matrix whose eigenvalues lie in $O_k^{\times}$. The norm $\n{\cdot}_I$ on $k^n$ associated to a basis $I \subset k^n$ giving a Jordan normal form of $U$ makes $U$ a unitary operator on the Banach $k$-vector space $(k^n,\n{\cdot}_I)$. Then $\n{\cdot}_I$ does not coincide with the original norm $\n{\cdot}$ and $U$ is not unitary with respect to $\n{\cdot}$, unless the matrix representing the change of basis of $k^n$ from the canonical basis to $I$ is not unitary.

\section{Functional Calculus and Unitary Diagonalisability}
\label{Functional Calculus and Unitary Diagonalisability}

Let $\A$ denote a Banach $k$-algebra in the following in this paper. We define a spectrum of an $A \in \A$ in \S \ref{The Spectrum of a Bounded Operator}, and formulate several notions of functional calculi of $A$ in $\A$ in \S \ref{Rigid Continuous Functional Calculus} and \S \ref{Reductive Functional Calculus}. Roughly speaking, a functional calculus of $A$ is a substitution of $A$ for a functions on its spectrum for which the coordinate function corresponds to $A$ itself. In order to ensure the uniqueness of such a substitution, we use the class of rigid continuous functions introduced in \S \ref{Rigidity of Continuous Functions}.

\subsection{The Spectrum of a Bounded Operator}
\label{The Spectrum of a Bounded Operator}

In order to define the normality of an operator, we formulate the notion of a functional calculus with respect to the class of rigid continuous functions. For this sake, we recall the notion of the spectrum of an element of a Banach $k$-algebra. The spectrum of an element of Banach $k$-algebra is an analogue of the spectrum of a bounded operator on an Archimedean Banach space, and naturally yields a notion of the spectrum of a bounded operator of a bounded operator on a Banach $k$-vector space.

\begin{dfn}
Let $R$ be a ring, and $S$ an $R$-algebra. For each $A \in S$, we set $\sigma_S(A) \coloneqq \set{\lambda \in R}{A - \lambda \in S \backslash S^{\times}}$, and call $\sigma_S(A)$ (resp.\ $k \backslash \sigma_S(A)$) {\it the spectrum} (resp.\ {\it the resolvent}) {\it of $A$ in $S$}. We abbreviate ``in $S$'' in the convention as long as there is no ambiguity.
\end{dfn}

We consider the case where $R$ is $k$ and $S$ is the underlying $k$-algebra of $\A$. For an $A \in \A$, {\it the spectrum $\sigma_{\A}(A) \subset k$ of $A$ in $\A$} is the spectrum of $A$ in the underlying $k$-algebra of $\A$.

\begin{dfn}
An $A \in \A$ is said to be {\it reductively finite} (resp.\ {\it reductively infinite}) if $\n{A} \in \v{k^{\times}}$ and if $\sigma_{\overline{\A}}(c^{-1}A + \A(1-))$ is a finite set (resp.\ an infinite set) for a $c \in k^{\times}$ with $\n{A} = \v{c}$.
\end{dfn}

The cardinality of $\sigma_{\overline{\A}}(c^{-1}A + \A(1-))$ is independent of the choice of $c \in k^{\times}$ with $\n{A} = \v{c}$, because a replacement of $c$ causes a replacement of $c^{-1}A + \A(1-)$ by an element in $\overline{k}{}^{\times} (c^{-1}A + \A(1-))$.

\begin{dfn}
Let $V$ be a Banach $k$-vector space. An operator $A$ on $V$ is said to be {\it reductively finite} (resp.\ {\it reductively infinite}) if $A$ is bounded and is reductively finite (resp.\ reductively infinite) in $\B_k(V)$.
\end{dfn}

The most important example of a reductively finite operator is a matrix. Suppose that $\A$ is the matrix algebra $\t{M}_n(k)$. Every $M \in \t{M}_n(k)$ is integral over $k$, but its characteristic polynomial is not necessarily decomposed in $k$ into linear functions when $k$ is not algebraically closed. Therefore $\sigma_{\t{M}_n(k)}(M)$ might be empty. In the infinite dimensional case, there is more complicated phenomenon peculiar to non-Archimedean analysis.

\begin{exm}
We consider the residue field $\overline{\C}_p$ of $\C_p$. Let $\lbrack \cdot \rbrack \colon \overline{\C}_p \to \C_p$ denote the Teichm\"uller embedding, which is unique because $\overline{\C}_p$ consists of $0$ and roots of unity. We set $I \coloneqq (\overline{\C}_p \sqcup \ens{\infty}) \times \N$. Then $I$ is a countable discrete topological space, and $\t{C}_0(I,\C_p)$ is a strictly Cartesian $\C_p$-Banach space of countable type with an orthonormal Schauder basis $\set{\delta_i}{i \in I}$. We consider the operator
\begin{eqnarray*}
  A \colon \C_p^{\oplus I} & \to & \t{C}_0(I,\C_p) \\
  \delta_i & \mapsto &
  \left\{
    \begin{array}{ll}
      \delta_{(\infty,n+1)} & (i = (\infty,n) \in \ens{\infty} \times \N) \\
      \delta_{(0,n-1)} & (i = (0,n) \in \ens{0} \times (\N \backslash \ens{0})) \\
      1 + \lbrack c \rbrack \delta_{(\overline{c},0)} & (i = (\overline{c},0) \in \overline{\C}_p \times \ens{0}) \\
      \delta_{(\overline{c},n-1)} + \lbrack \overline{c} \rbrack \delta_{(\overline{c},n)} & (i = (\overline{c},n) \in \overline{\C}_p^{\times} \times (\N \backslash \ens{0}))
    \end{array}
  \right.
\end{eqnarray*}
densely defined on $\t{C}_0(I,\C_p)$. An estimation of the norms of entries of the matrix representation, $A$ uniquely extends to an everywhere defined bounded operator on $\t{C}_0(I,\C_p)$ satisfying $\sigma_{\B_{\C_p}(\t{C}_0(I,\C_p))}(A) \subset O_{\C_p}$. Let $c \in O_{\C_p}$. We put $\overline{c} \coloneqq c + O_{\C_p} \in \overline{C}_p$ and $\varpi_c \coloneqq c - \lbrack c + O_{\C_p} \rbrack \in O_{\C_p}$. For each $(n,m,c_0) \in \N \times \N \times \N \times O_{\C_p}$ with $n-1 \geq m$, we set
\begin{eqnarray*}
  \mu_{c_0}^{n,m} \coloneqq \sum_{l = 0}^{n-m-1} \frac{l!}{(l-m)! m!} (\lbrack \overline{c} \rbrack - c_0)^{n-l-2}(-c_0)^{l-m} \in O_{\C_p}.
\end{eqnarray*}
We consider the operator
\begin{eqnarray*}
  B_c \colon \C_p^{\oplus I} & \to & \t{C}_0(I,\C_p) \\
  \delta_i & \mapsto &
  \left\{
    \begin{array}{ll}
      \sum_{m = 0}^{n-1} c^{n-m-1} \delta_{(\infty,m)} + c^n \sum_{m = 0}^{\infty} \varpi_c^m \delta_{(\overline{c},m)} & (i = (\infty,n) \in \ens{\infty} \times \N) \\
      \sum_{m = n+1}^{\infty} \varpi_c^{m-n-1} \delta_{(\lbrack \overline{c} \rbrack,m)} & (i = (\overline{c},n) \in \ens{\overline{c}} \times \N) \\
      (\lbrack \overline{c} \rbrack - \lbrack \overline{c}_0 \rbrack)^{-n} \delta_{(\lbrack \overline{c} \rbrack,0)} - \sum_{m = 0}^{n-1} \mu_{\lbrack \overline{c}_0 \rbrack}^{n,m} \delta_{(\overline{c}_0,m)} & (i = (\overline{c}_0,n) \in (\overline{\C}_p \backslash \ens{\overline{c}}) \times \N)
    \end{array}
  \right.
\end{eqnarray*}
densely defined on $\t{C}_0(I,\C_p)$, which uniquely extends to an everywhere defined bounded operator on $\t{C}_0(I,\C_p)$ by a similar argument. Applying $A$ and $B_c$ to $\delta_i$ for each $i \in I$, we obtain $(A-c)B_c = B_c(A-c) = 1 \in \B_{\C_p}(\t{C}_0(I,\C_p))$, and hence $A-c \in \B_{\C_p}(\t{C}_0(I,\C_p))^{\times}$. It implies $\sigma_{\B_{\C_p}(\t{C}_0(I,\C_p))}(A) = \emptyset$.
\end{exm}

Such a phenomenon occurs because Gel'fand--Mazur theorem, which holds for $\C$, does never hold for a valuation field. We regard an operator with empty spectrum as exceptional in this paper, because such an operator never admits a functional calculus. In order to give a criterion for an operator with enough points in spectra, we study the reduction. Before that, we verify several basic properties of the spectrum.

\begin{prp}
\label{functoriality 2}
Let $R$ be a ring.
\begin{itemize}
\item[(i)] Let $S_1$ and $S_2$ be $R$-algebras. For any $R$-algebra homomorphism $\phi \colon S_1 \to S_2$ and $A \in S_1$, $\sigma_{S_2}(\phi(A))$ is contained in $\sigma_{S_1}(A)$.
\item[(ii)] Let $R_0 \subset R$ be a subring, and $S$ an $R$-algebra. For any $A \in S$, $\sigma_{S_0}(A)$ coincides with $\sigma_S(A) \cap R_0$, where $S_0$ denotes $S$ regarded as an $R_0$-algebra.
\end{itemize}
\end{prp}

\begin{proof}
The first assertion follows from $\phi(S_1^{\times}) \subset S_2^{\times}$, and the second assertion follows from the definition of the spectrum.
\end{proof}

\begin{lmm}
\label{inverse}
Let $\A$ be a Banach $k$-algebra. Then $\A(1)^{\times}$ and $\A^{\times}$ are open subsets of $\A$ containing $1 + \A(1-)$.
\end{lmm}

\begin{proof}
The convergent radius of $(1 - T)^{-1} = \sum_{i = 0}^{\infty} T^i \in \A(1)[[T]]$ is $1-$, and hence $1 + \A(1-)$ is contained in $\A(1)^{\times}$. Put $G \coloneqq \A(1)^{\times}$ (resp.\ $G \coloneqq \A^{\times}$). Let $f \in G$. Then for any $f' \in \A(1)$ (resp.\ $f' \in \A$) with $\n{f-f'} < \n{f^{-1}}^{-1}$, we have $f' = f + (f'-f) = f(1 + f^{-1}(f'-f)) \in G$, and hence $G$ is an open subset of $\A$.
\end{proof}

\begin{prp}
\label{spectrum -> bounded closed}
For any $A \in \A$, $\sigma_{\A}(A)$ is a closed subset of $k$ contained in $k(\n{A})$.
\end{prp}

\begin{proof}
For any $c \in k$ with $\v{c} > \n{A}$, we have $\n{a^{-1}A} < 1$ and hence $A-c = -c(1-c^{-1}A) \in k^{\times}(1 + \A(1-)) \subset \A^{\times}$ by Lemma \ref{inverse}. Therefore $\sigma_{\A}(A)$ is contained in $k(\n{A})$. Let $c \in k \backslash \sigma_{\A}(A)$. For any $c' \in k$ with $\v{c'-c} < \n{(A-c)^{-1}}$, we have $\n{(c'-c)(A-c)^{-1}} < 1$, and hence $A-c' = (A-c)-(c'-c) = (A-c)(1-(c'-c)(A-c)^{-1}) \in \A^{\times}(1 + \A(1-)) \subset \A^{\times}$ by Lemma \ref{inverse}. Therefore $\sigma_{\A}(A)$ is closed in $k$.
\end{proof}

\begin{crl}
\label{spectrum -> bounded closed  2}
Suppose that $k$ is a finite field endowed with the trivial norm or a local field. For any $A \in \A$, $\sigma_{\A}(A)$ is a compact subspace of $k$.
\end{crl}

\begin{proof}
The assertion follows from Proposition \ref{spectrum -> bounded closed} because every bounded closed subset of a finite field or a local field is compact.
\end{proof}

\begin{crl}
\label{spectrum -> integral spectrum}
For any $A \in \A(1)$, $\sigma_{\A}(A)$ is contained in $\sigma_{\A(1)}(A)$.
\end{crl}

\begin{proof}
Let $\A_0$ denotes $\A$ regarded as an $O_k$-algebra. We have $\sigma_{\A}(A) \subset O_k$ by Proposition \ref{spectrum -> bounded closed}, $\sigma_{\A_0}(A) = \sigma_{\A}(A) \cap O_k = \sigma_{\A}(A)$ by Proposition \ref{functoriality} (ii), and $\sigma_{\A}(A) = \sigma_{\A_0}(A) \subset \sigma_{\A(1)}(A)$ by Proposition \ref{functoriality 2} (i). 
\end{proof}

We emphasise that $\sigma_{\A}(A)$ depends not only on $A$ but also $\A$, while the spectrum of an element of an Archimedean $C^*$-algebra does not depend on the algebra. The proof of the independence of the choice of $\A$ in the Archimedean case heavily uses properties peculiar to the topologies of $\R$ and $\C$. See \cite{Dou72} 4.28 for more detail.

\begin{exm}
\label{spectrum of affinoids}
The equalities $\sigma_{k \ens{T}}(T) = O_k$ and $\sigma_{k \ens{T,T^{-1}}}(T) = O_k^{\times}$ hold, and hence the spectrum of $T$ depends on the choice of a Banach $k$-algebra to which $T$ belongs.
\end{exm}

\subsection{Rigid Continuous Functional Calculus}
\label{Rigid Continuous Functional Calculus}

We introduce the notion of a functional calculus with respect to the class of rigid continuous function. As is mentioned in \S \ref{Introduction}, a functional calculus means a substitution of a substitution of an operator to a function on its spectrum. In other words, a functional calculus gives a way to regard an operator as a function on a topological space. Let $A$ denote an element of $\A$ in the following.

\begin{dfn}
\label{rcfc}
{\it The rigid continuous functional calculus of $A$ in $\A$} is an isometric $k$-algebra homomorphism $\iota_A \colon \t{C}_{\t{rig}}(\sigma_{\A}(A),k) \to \A$ sending $z_{\sigma}$ to $A$.
\end{dfn}

The rigid continuous functional calculus of $A$ in $\A$ is unique because the localisation of $k[z] \subset \t{C}_{\t{rig}}(\sigma_{\A}(A),k)$ is dense. The rigid continuous functional calculus is deeply related to unitary operators through the condition of isometry. When the equality $\t{C}_{\t{rig}}(\sigma_{\A}(A),k) = \t{C}_{\t{bd}}(\sigma_{\A}(A),k)$ holds, then we simply call $\iota_A$ {\it the continuous functional calculus of $A$}.

\begin{rmk}
\label{non-empty spectrum}
If $A$ admits the rigid continuous functional calculus in $\A$, then the equality $\t{C}_{\t{rig}}(\emptyset,k) = \ens{0}$ ensures that $\sigma_{\A}(A)$ is non-empty unless $\A = \ens{0}$. The image of $\t{C}_{\t{rig}}(\sigma_{\A}(A),k)$ by $\iota_A$ is a closed $k$-subalgebra of $\Leb_{\A}(A)$ because $\t{C}_{\t{rig}}(\sigma_{\A}(A),k)$ is complete and $\iota_A$ is an isometry. If $\sigma_{\A}(A)$ is compact, then we have $\t{C}_{\t{rig}}(\sigma_{\A}(A),k) = \t{C}(\sigma_{\A}(A),k)$ by Proposition \ref{continuous functions on a compact subset}, and the continuous functional calculus induces an isometric isomorphism $\t{C}_{\t{rig}}(\sigma_{\A}(A),k) \cong \Leb(A)$.
\end{rmk}

\begin{prp}
\label{spectral}
If $A$ admits the rigid continuous functional calculus in $\A$, then the equality $\n{A} = \sup_{\lambda \in \sigma_{\A}(A)} \v{\lambda}$ holds.
\end{prp}

\begin{proof}
The right hand side coincides with the supremum norm of $z_{\sigma_{\A}(A)}$, and hence the equality holds because $\iota_A$ is an isometry.
\end{proof}

\begin{dfn}
\label{normal}
A bounded operator on a Banach $k$-vector space $V$ is said to be {\it normal} if it admits the rigid continuous functional calculus in $\B_k(V)$.
\end{dfn}

We note that another formulation of the normality of is dealt with in \cite{Koc13} \S 2, where a normal operator is assumed to admit the continuous functional calculus. In particular, a normal operator in the sense of \cite{Koc13} is a normal operator in our sense.

\begin{prp}
\label{n <-> ud}
Let $n \in \N$ An $M \in \t{M}_n(k)$ is normal if and only if $M$ is diagonalisable by a unitary matrix in $k$.
\end{prp}

\begin{proof}
Put $\sigma \coloneqq \sigma_{\t{M}_n(k)}(M)$. Since $\t{M}_n(k)$ is of finite dimension, $\sigma$ is a finite set, and hence $\t{C}_{\t{rig}}(\sigma,k)$ coincides with $\t{C}(\sigma,k)$ by Proposition \ref{continuous functions on a compact subset}. Suppose that $M$ is a normal operator with the rigid continuous functional calculus $\iota_M$. For each $\lambda \in \sigma$, we denote by $V_{\lambda} \subset k^n$ the range of $\iota_M(\delta_{\lambda})$. Since $(\delta_{\lambda})_{\lambda \in \sigma}$ is a partition of unity in $\t{C}(\sigma,k)$, we obtain an orthogonal decomposition $V = \bigoplus_{\lambda \in \sigma} V_{\sigma}$ by Proposition \ref{orthogonal property}. Since $\iota_M(\delta_{\lambda})$ commutes with $\iota_M(z_{\sigma}) = M$, the decomposition $V = \bigoplus_{\lambda \in \sigma} V_{\sigma}$ is compatible with the action of $M$, and the restriction of $M$ to $V_{\lambda}$ coincides with the scalar multiplication by $\lambda$ for any $\lambda \in \sigma$. By \cite{BGR84} 2.5.1/5, every subspace of $k^n$ is a strictly Cartesian Banach $k$-vector space and admits an orthonormal basis. Therefore we obtain an orthonormal basis $I$ of $k^n$ given as the collection of orthonormal bases of $V_{\lambda}$ for each $\lambda \in \sigma$. By Proposition \ref{criterion for unitarity}, the matrix $U$ whose columns are given as elements of $I$ is a unitary matrix, and $U^{-1} M U$ is a diagonal matrix.

\vspace{0.2in}
Suppose that $M$ is diagonalisable by a unitary matrix $U \in \t{M}_n(k)$. For each $\lambda \in \sigma$, let $V_{\lambda} \subset k^n$ denote the eigenspace of $M$ belonging to the eigenvalue $\lambda$, and $P_{\lambda} \colon k^n \twoheadrightarrow V_{\lambda}$ the projection induced by the decomposition $V = \bigoplus_{\lambda' \in \sigma} V_{\lambda'}$. Since $U^{-1} M U$ is a diagonal matrix, the collection $I$ of columns of $U$ is an orthonormal basis consisting of eigenvectors of $M$ by Proposition \ref{criterion for unitarity}. Therefore we obtain $\n{v} = \max_{\lambda \in \sigma} \n{P_{\lambda} v}$ for any $v \in k^n$. The inequality $V_{\lambda} \neq 0$ ensures $\n{P_{\lambda}} = 1$ for any $\lambda \in \sigma$, and hence $(P_{\lambda})_{\lambda \in \sigma}$ forms an orthonormal partition of unity in $\t{M}_n(k)$. We define a map $\iota_M \colon k^{\sigma} \to \t{M}_n(k)$ by setting $\iota_M((c_{\lambda})_{\lambda \in \sigma}) \coloneqq \sum_{\lambda \in \sigma} \lambda P_{\lambda}$ for each $(c_{\lambda})_{\lambda \in \sigma} \in k^{\sigma}$. Since $(P_{\lambda})_{\lambda \in \sigma}$ is orthonormal, $\iota_M$ is an isometry. The equality $M = \sum_{\lambda \in \sigma} M P_{\lambda} = \sum_{\lambda \in \sigma} \lambda P_{\lambda}$ ensures $\iota_M(z_{\sigma}) = M$, and hence $\iota_M$ is the continuous functional calculus of $M$ in $\t{M}_n(k)$.
\end{proof}

\begin{rmk}
Let $n \in \N$. As we mentioned in the end of \S \ref{Operator Algebra}, the notion of the unitarity of an operator on $k^n$ makes sense only when we fix a norm of $k^n$. Therefore a replacement of the canonical basis of $k$ by another basis does not necessarily preserve the normality of an operator. Moreover, a matrix is diagonalisable in $k$ if and only if it represents a normal bounded operator on $k^n$ with respect to some norm of $k^n$.
\end{rmk}

Imitating \cite{Ber90} 7.3.2 Theorem, we define another functional calculus with stronger rigidity than the rigid continuous functional calculus. We denote by $\Aff^1_k$ the rigid analytic affine line over $k$, which is canonically isomorphic to the analytification of $\t{Spec}(k[T])$. For conventions for Berkovich spaces, see \cite{Ber90}. For each affinoid domain $V \subset \Aff^1_k$, we regard the underlying $k$-vector space $\t{H}^0(V,\Opn_{\Aff^1_k})$ as a $k$-subalgebra of $\t{C}_{\t{bd}}(V(k),k)$ through the submetric injective $k$-algebra homomorphism $\t{H}^0(V,\Opn_{\Aff^1_k}) \hookrightarrow \t{C}_{\t{bd}}(V(k),k)$ given by the evaluation at each $k$-rational point. In particular, for a rational domain of a closed disc in $\Aff^1_k$, we regard the underlying $k$-algebra of $\t{H}^0(V,\Opn_{\Aff^1_k})$ as a $k$-subalgebra of $\t{C}_{\t{rig}}(V(k),k)$. For a closed subset $\sigma \subset k$, a {\it rational neighbourhood} of $\sigma$ is a rational domain of a closed disc in $\Aff^1_k$ containing the image of $\sigma$ by the identification $k \cong \Aff^1_k(k)$.

\begin{dfn}
For a rational neighbourhood $V \subset \Aff^1_k$ of $\sigma_{\A}(A)$, {\it the holomorphic functional calculus of $A$ in $\A$ on $V$} is a bounded $k$-algebra homomorphism $\iota_A \colon \t{H}^0(V, \Opn_{\Aff^1_k}) \to \A$ sending $z_{V(k)}$ to $A$. A {\it holomorphic functional calculus of $A$ in $\A$} is the holomorphic functional calculus of $A$ in $\A$ on $V$ for some rational neighbourhood $V \subset \Aff^1_k$ of $\sigma_{\A}(A)$.
\end{dfn}

For a fixed rational neighbourhood $V$ of $\sigma_{\A}(A)$, the holomorphic functional calculus of $A$ in $\A$ on $V$ is unique because the localisation of $k[T]$ is dense in the corresponding affinoid domain. The restriction of a rigid analytic function on an affinoid domain to a smaller affinoid domain is again rigid analytic, and hence the rational neighbourhood appearing in the definition of a holomorphic functional calculus can be replaced by a larger one.

\begin{exm}
Let $A \in \A$ be an element admitting the rigid continuous functional calculus. For any rational neighbourhood $V \subset \Aff^1_k$ of $\sigma_{\A}(A)$, the composite of the evaluation map $\t{H}^0(V,\Opn_{\Aff^1_k}) \to \t{C}_{\t{rig}}(V(k),k)$, the restriction map $\t{C}_{\t{rig}}(V(k),k) \to \t{C}_{\t{rig}}(\sigma_{\A}(A),k)$, and $\iota_A$ sends $z_{V(k)}$ to $A$, and hence is the holomorphic functional calculus of $A$ in $\A$ on $V$.
\end{exm}

\begin{exm}
Suppose that the valuation of $k$ is non-trivial. Let $V_0 \subset \Aff^1_k$ be an affinoid domain, and $V \subset \Aff^1_k$ a rational domain of a closed disc in $\Aff^1_k$ with $V_0 \subset V$. Then $V$ is a rational neighbourhood of $V_0(k)$, and $V_0$ itself is a rational domain of a closed disc in $\Aff^1_k$ by \cite{Ber90} 2.6.4 Remark. Therefore $V_0(k)$ coincides with $\sigma_{\t{H}^0(V_0, \Opn_{\Aff^1_k})}(z_{V(k)})$ by \cite{Ber90} 7.1.5 Corollary. The restriction map $\t{H}^0(V, \Opn_{\Aff^1_k}) \to \t{H}^0(V_0, \Opn_{\Aff^1_k})$ sends $z_{V(k)}$ to $z_{V_0(k)}$, and hence is the holomorphic functional calculus of $z_{V_0(k)}$ in $\t{H}^0(V_0, \Opn_{\Aff^1_k})$ on $V$.
\end{exm}

\subsection{Reductive Functional Calculus}
\label{Reductive Functional Calculus}

In order to establish an algorithm for a criterion of the unitarity of an operator, we study the functional calculus of the reduction.

\begin{dfn}
An $A \in \A$ is said to {\it admit the reductive functional calculus in $\A$} if $\n{A} \in \v{k^{\times}}$ and if for a $c \in k^{\times}$ with $\n{A} = \v{c}$, $c^{-1}A + \A(1-)$ admits the continuous functional calculus in $\overline{\A}$.
\end{dfn}

Let $A \in \A$ be an element admitting the reductive functional calculus. Then for any $c \in k^{\times}$ with $\n{A} = \v{c}$, $c^{-1}A + \A(1-)$ admits the continuous functional calculus, because a replacement of $c$ causes a replacement of $c^{-1}A + \A(1-)$ by an element in $\overline{k}{}^{\times}(c^{-1}A + \A(1-))$.

\begin{dfn}
Let $V$ be a Banach $k$-vector space. A bounded operator $A$ on $V$ is said to {\it reductively normal} if $\n{A} \in \v{k^{\times}}$ and if for a $c \in k^{\times}$ with $\n{A} = \v{c}$, $c^{-1}A$ induces a normal operator on $\overline{V}$ through the canonical projection $V \twoheadrightarrow \overline{V}$.
\end{dfn}

\begin{prp}
\label{rn <-> rfc}
Let $V$ be a strictly Banach $k$-vector space. Suppose that $V$ is of finite dimension or that $k$ is a finite field endowed with the trivial valuation or a local field. A bounded operator on $V$ is reductively normal if and only if it admits the reductive functional calculus.
\end{prp}

\begin{proof}
If $V$ is of finite dimension, then so is $\overline{V}$, and hence the spectrum of an everywhere-defined operator on $\overline{V}$ is a finite set. If $k$ is a finite field endowed with the trivial valuation or a local field, then $\overline{k}$ is a finite field endowed with the trivial valuation, and hence the spectrum of an everywhere-defined operator on $\overline{V}$ is a finite set. Therefore the spectrum of an everywhere-defined operator on $\overline{V}$ is a finite set in both cases, and the rigid continuous functional calculus of an everywhere-defined operator on $\overline{V}$ is always the continuous functional calculus. Thus the assertion follows from Proposition \ref{reduction & operator reduction}.
\end{proof}

\begin{exm}
Let $n \in \N$. An $M \in \t{M}_n(k)$ admits the reductive functional calculus in $\t{M}_n(k)$ if and only if $M \neq 0$ and for a $c \in k^{\times}$ with $\n{M} = \v{c}$, $\Pi_{k^n}(c^{-1}M) \in \t{M}_n(\overline{k})$ is diagonalisable in $\overline{k}$.
\end{exm}

We study the relation between functional calculi and reductions. For this, we show a certain compatibility of spectra and reductions. For an $A \in \A(1)$, the image of the spectrum of an element $A$ by the reduction $O_k \twoheadrightarrow \overline{k}$ does not necessarily coincide with the spectrum of the image of $A$ by the reduction $\A \twoheadrightarrow \overline{\A}$. On the other hand, the following ensure that the spectrum of $A$ in $\A(1)$ is sufficiently large.

\begin{lmm}
\label{inverse and reduction}
The group homomorphism $\A(1)^{\times} \to \overline{\A}{}^{\times}$ induced by the canonical projection $\A(1) \twoheadrightarrow \overline{\A}$ is surjective, and an $f \in \A(1)$ lies in $\A(1)^{\times}$ if and only if $f + A(1-) \in \overline{\A}$ lies in $\overline{\A}{}^{\times}$.
\end{lmm}

\begin{proof}
It suffices to verify that every $f \in \A(1)$ with $f + A(1-) \in A + \A(1-){}^{\times}$ lies in $\A(1-)^{\times}$. Let $f \in \A(1)$ with $f + A(1-) \in A + \A(1-){}^{\times}$. Take a lift $g \in \A(1)$ of $(f + A(1-))^{-1}$. We have $(fg - 1) + A(1-) = (f + A(1-))(f + A(1-))^{-1} - (1 + A(1-)) = 0 \in A + \A(1-)$, and hence $\n{ab - 1} < 1$. We obtain $fg \in 1 + \A(1-) \subset \A(1)^{\times}$ by Lemma \ref{inverse}. It implies $f \in \A(1)^{\times}$.
\end{proof}

\begin{prp}
\label{spectrum and reduction}
For any $A \in \A(1)$, the equality $\sigma_{\A(1)}(A) = \bigcup_{\Lambda \in \sigma_{\overline{\A}}(A + \A(1-))} \Lambda$ holds.
\end{prp}

\begin{proof}
The assertion immediately follows from Lemma \ref{inverse and reduction}.
\end{proof}

\begin{crl}
\label{spectrum and reduction 2}
For any $A \in \A(1)$, $\sigma_{\overline{\A}}(A + \A(1-))$ coincides with the image of $\sigma_{\A(1)}(A)$ by the reduction $O_k \twoheadrightarrow \overline{k}$.
\end{crl}

\begin{crl}
\label{spectrum and reduction 3}
For any non-empty closed subset $\sigma \subset k$ contained in $O_k$, $\sigma_{\overline{\t{C}_{\t{rig}}(\sigma,k)}}(z^{(k)}_{\sigma} + \t{C}_{\t{rig}}(\sigma,k)(1-))$ coincides with the image of $\sigma$ by the reduction $O_k \twoheadrightarrow \overline{k}$.
\end{crl}

\begin{proof}
The assertion immediately follows from Corollary \ref{invertible rigid continuous function 2} and Corollary \ref{spectrum and reduction 2}.
\end{proof}

\begin{crl}
\label{spectrum and reduction 4}
For any $A \in \A(1)$ admitting the rigid continuous functional calculus in $\A$, $\sigma_{\overline{\A}}(A + \A(1-))$ coincides with the image of $\sigma_{\A}(A)$ by the reduction $O_k \twoheadrightarrow \overline{k}$.
\end{crl}

\begin{proof}
We have $\sigma_{\A}(A) \subset \sigma_{\A(1)}(A) \subset \set{\lambda + \varpi}{(\lambda,\varpi) \in \sigma_{\A}(A) \times m_k}$ by Proposition \ref{functoriality 2} and Corollary \ref{invertible rigid continuous function 2}. Therefore the assertion immediately follows from Corollary \ref{spectrum and reduction 3}.
\end{proof}

As a result, we obtain a necessary condition for the existence of the rigid continuous functional calculus using an information of the reduction.

\begin{lmm}
\label{n implies rn}
Let $A \in \A(1)$. If $A$ admits the rigid continuous functional calculus in $\A$, then $A + \A(1-)$ admits the rigid continuous functional calculus in $\overline{\A}$.
\end{lmm}

\begin{proof}
Put $\sigma \coloneqq \sigma_{\A}(A) \subset O_k$. We denote by $\overline{\sigma} \subset \overline{k}$ the image of $\sigma$ by the reduction $O_k \twoheadrightarrow \overline{k}$. We have $\sigma_{\overline{\A}}(A + \A(1-)) = \overline{\sigma}$ by Corollary \ref{spectrum and reduction 2}. Let $\iota_A \colon \t{C}_{\t{rig}}(\sigma,k) \hookrightarrow \A$ denote the rigid continuous functional calculus of $A$ in $\A$. Let $\overline{f} \in \t{C}_{\t{rig}}(\overline{\sigma},\overline{k})$. Since $\overline{k}$ is discrete, so is $\t{C}_{\t{bd}}(\overline{\sigma},\overline{k})$, and hence $\t{C}_{\t{rig}}(\overline{\sigma},\overline{k})$ coincides with the $\overline{k}$-subalgebra of $\t{C}(\overline{\sigma},\overline{k})$ given by the localisation of $\Leb(z^{(\overline{k})}_{\overline{sigma}}) = \overline{k}[z^{(\overline{k})}_{\overline{sigma}}]$ by $\overline{k}[z^{(\overline{k})}_{\overline{sigma}}] \cap \t{C}(\overline{\sigma},\overline{k})^{\times}$. Therefore there is an $(\overline{F}_0, \overline{F}_{\infty}) \in \overline{k}[T] \times \overline{k}[T]$ with $\overline{F}_{\infty}(z^{(\overline{k})}_{\overline{sigma}}) \in \overline{k}[z^{(\overline{k})}_{\overline{sigma}}] \cap \t{C}(\overline{\sigma},\overline{k})^{\times}$ and $\overline{F}_{\infty}(z^{(\overline{k})}_{\overline{sigma}})^{-1} \overline{F}_0(z^{(\overline{k})}_{\overline{sigma}}) = \overline{f}$. Take a lift $(F_0,F_{\infty}) \in O_k[T] \times O_k[T]$ of $(\overline{F}_0, \overline{F}_{\infty})$. Since $\overline{F}_{\infty}(z^{(\overline{k})}_{\overline{sigma}})$ is invertible in $\t{C}(\overline{\sigma},\overline{k})$, it has no zero in $\overline{\sigma}$. Therefore $F_{\infty}(z^{(k)}_{\sigma})$ factors through the inclusion $O_k^{\times} \hookrightarrow k$ by Corollary \ref{spectrum and reduction 3}. It implies that $F_{\infty}(z^{(k)}_{\sigma})$ is invertible in $\t{C}_{\t{bd}}(\sigma,k)(1)$, and $F_{\infty}(z^{(k)}_{\sigma})^{-1}$ lies in $\t{C}_{\t{rig}}(\sigma,k)(1)^{\times}$. Put $F(z^{(k)}_{\sigma}) \coloneqq F_{\infty}(z^{(k)}_{\sigma})^{-1} F_0(z^{(k)}_{\sigma}) \in \t{C}_{\t{rig}}(\sigma,k)(1)$. Then the composite of $F(z^{(k)}_{\sigma})$ and the reduction $O_k \twoheadrightarrow \overline{k}$ coincides with the composite of the canonical projection $\sigma \twoheadrightarrow \overline{\sigma}$ and $\overline{f}$, and hence $\overline{f}(A + \A(1-)) \coloneqq \iota_A(F(z^{(k)}_{\sigma})) + \A(1-) \in \overline{\A}$ depends only on $\overline{f}$. Since $\iota_A$ is an isometry, we have $\overline{f}(A + \A(1-)) \neq 0$ unless $\overline{f} = 0$ by Corollary \ref{spectrum and reduction 3}. Therefore the correspondence $\overline{f} \mapsto \overline{f}(A + \A(1-))$ gives an injective $\overline{k}$-algebra homomorphism $\iota_{A + \A(1-)} \colon \t{C}_{\t{rig}}(\overline{\sigma},\overline{k}) \hookrightarrow \overline{\A}$ sending $z^{(\overline{k})}_{\overline{\sigma}}$ to $A + \A(1-)$. Thus $A + \A(1-)$ admits the rigid continuous functional calculus in $\overline{\A}$.
\end{proof}

\begin{thm}
\label{n implies rn 2}
Let $A \in \A \backslash \ens{0}$. If $\sigma_{\A}(A)$ is compact and $A$ admits the continuous functional calculus in $\A$, then $A$ admits the reductive functional calculus in $\A$.
\end{thm}

\begin{proof}
There is a $c \in \sigma_{\A}(A)$ with $\v{c} = \sup_{\lambda \in \sigma_{\A}(A)} \v{\lambda}$ by the maximal modulus principle. Since $A$ admits the continuous functional calculus in $\A$, we have $\n{A} = \n{z^{(k)}_{\sigma_{\A}(A)}} = \v{c}$, and hence $c^{-1}A$ is an element with $\n{c^{-1}A} = 1$ admitting the continuous functional calculus in $\A$. Since $\sigma_{\A}(A)$ is compact, so is $\sigma_{\A}(c^{-1}A) = \set{c^{-1}\lambda}{\lambda \in \sigma_{\A}(A)}$, and so is $\sigma_{\overline{\A}}(c^{-1}A + \A(1-))$ by Corollary \ref{spectrum and reduction 2}. Therefore the assertion follows from Proposition \ref{continuous functions on a compact subset} and Lemma \ref{n implies rn}.
\end{proof}

\begin{crl}
\label{n implies rn 3}
Let $V$ be a strictly Cartesian Banach $k$-vector space. Suppose that $V$ is of finite dimension or that $k$ is a finite field endowed with the trivial norm or a local field. Then every non-zero normal bounded operator on $V$ is reductively normal.
\end{crl}

\begin{proof}
The assertion follows from Corollary \ref{spectrum -> bounded closed  2}, Proposition \ref{rn <-> rfc}, and Theorem \ref{n implies rn 2}.
\end{proof}

\section{Criterion of the Normality}
\label{Criterion of the Normality}

In this section, we study the relation between the normality and the operator reduction. A bounded operator on a Banach $k$-vector space is not necessarily normal even if its operator reduction is normal, but we establish a criterion of the normality including ``the repetitive reduction'' in \S \ref{Reductively Finite Operator}. Although such a criterion contains an infinitely many recursive computation, it yields an algorithm when the representation space is of finite dimension. Indeed, we give an explicit algorithm for a criterion of the diagonalisability of a matrix by a unitary matrix in \S \ref{Reduction Algorithm for Matrices}.

\subsection{Reductively Infinite Operator}
\label{Reductively Infinite Operator}

We deal with a reductively infinite operator in this subsection. The spectrum of a reductively infinite operator is mostly non-compact, and hence Vishik's non-Archimedean analogue of Riesz functional calculus in \cite{Vis85} does not work here. When $\overline{k}$ is a finite field, then every bounded operator is reductively finite. Therefore we assume that $\overline{k}$ is an infinite field throughout this subsection.

\begin{exm}
Let $K/k$ be a non-trivial extension of complete valuation fields. Take an $r \in (0,1)$ with $K(r) \subsetneq m_K$, and put $\sigma \coloneqq O_K \backslash (O_k + K(r))$. We regard $\t{C}_{\t{bd}}(\sigma,K)$ as a Banach $k$-algebra through the embedding $k \hookrightarrow K$. Then we have $\sigma_{\t{C}_{\t{bd}}(\sigma,K)}(z_{\sigma}^{(K)}) = \emptyset$ and $\sigma_{\overline{\t{C}_{\t{bd}}(\sigma,K)}}(z_{\sigma}^{(K)} + \t{C}_{\t{bd}}(\sigma,K)(1-)) = \overline{k}$. Thus $z_{\sigma}^{(K)}$ is a reductively infinite element with compact spectrum.
\end{exm}

\begin{thm}
\label{hfc for a reductively infinite operator}
Let $A \in \A$ be a reductively infinite element. Then the $k$-algebra homomorphism $k[T] \to \A \colon T \mapsto A$ induces an isometric isomorphism $k \ens{\n{A}^{-1} T} \to \Leb(A)$ sending $T$ to $A$.
\end{thm}

\begin{proof}
Since $A$ is reductively infinite, there is a $c \in k^{\times}$ with $\n{A} = \v{c}$. Replacing $A$ by $c^{-1}A$, we may assume $\n{A} = 1$. We prove that $\n{F(A)} = 1$ for any $F \in O_k[T]$ with $\n{F} = 1$ with respect to the restriction of the Gauss norm of radius $1$ through the canonical embedding $k[T] \hookrightarrow k \ens{T}$. Since $F(A)$ is an $O_k$-linear combination of powers of $A$, we have $\n{F(A)} \leq 1$. We denote by $\overline{F} \in \overline{k}[T]$ the image of $F + A(1-) \in \overline{k \ens{T}}$ by the $\overline{k}$-algebra isomorphism in Example \ref{Tate algebra}. The equality $\n{F} = 1$ ensures $\overline{F} \neq 0 \in \overline{k}[T]$. Since $\overline{F}$ has at most finitely many zeros and $\sigma_{\overline{\A}}(A + A(1-))$ is an infinite set, there is a $\overline{\lambda} \in \sigma_{\overline{\A}}(A + A(1-))$ with $\overline{F}(\overline{\lambda}) \neq 0 \in \overline{k}$. We show $\overline{F}(A + \A(1-)) \neq 0 \in \overline{\A}$. Assume $\overline{F}(A + \A(1-)) = 0$. Let $\overline{G} \in \overline{k}[T]$ denote the unique element with $\overline{F}(T) - \overline{F}(\lambda) = (T - \lambda) \overline{G}$. We have $1 = - \overline{F}(\lambda)^{-1}(\overline{F}(A + \A(1-)) - \overline{F}(\lambda)) = - \overline{F}(\lambda)^{-1} \overline{G}(A + \A(1-))((A + \A(1-)) - \lambda)$, and hence $(A + \A(1-)) - \lambda \in \overline{\A}{}^{\times}$. It contradicts $\lambda \in \sigma_{\overline{\A}}(A + \A(1-))$. Therefore we obtain $\overline{F}(A + \A(1-)) \neq 0 \in \overline{\A}$. It implies $F(A) \in A(1-)$ and hence $\n{F(A)} = 1$.

\vspace{0.2in}
We verify that the $k$-algebra homomorphism $i_A \colon k[T] \to \A \colon T \mapsto A$ is an isometry. For this sake, it suffices to show that every $F \in k[T]$ with $\n{F(A)} \leq 1$ lies in $O_k[T]$. Take an $F \in k[T]$ with $\n{F(A)} \leq 1$. Assume $F \notin O_k[T]$. Then there is a $c \in k \backslash O_k$ such that $\n{F} = \v{c}$. We have $\n{c^{-1}F} = \v{c}^{-1} \n{F} = 1$ and hence $c^{-1}F \in O_k[T]$, which ensures $\n{(c^{-1}F)(A)} = 1$ by the argument above. It contradicts the inequality $\n{(c^{-1}F)(A)} = \v{c^{-1}} \ \n{F(A)} \leq \v{c}^{-1} < 1$. We obtain $F \in O_k[T]$. It implies that $i_A$ is an isometry. Since the image of $k[T]$ is dense in $k \ens{T}$, $i_A$ uniquely extends to an isometry $i_A \colon k \ens{T} \hookrightarrow \A$. We have $i_A(T) = A$ by definition. Since the image of $k[T]$ is dense in $k \ens{T}$, we obtain $i_A(k \ens{T}) \subset \Leb(A)$. Moreover, since $k[A]$ is dense in $\Leb(A)$, $i_A(k \ens{T})$ is dense in $\Leb(A)$. Since $k \ens{T}$ is complete and $i_A$ is an isometry, $i_A(k \ens{T})$ is closed in $k \ens{T}$. Therefore $i_A$ is onto $\Leb(A)$. We conclude that $i_A$ gives an isometric isomorphism $i_A \colon k \ens{T} \to \Leb(A)$.
\end{proof}

For each $r \in (0,\infty)$, we denote by $r \t{D}^1_k$ the closed disc in $\Aff^1_k$ of radius $r$.

\begin{crl}
Every reductively infinite element $A \in \A$ admits the holomorphic functional calculus in $\A$ on $\n{A} \t{D}^1_k$.
\end{crl}

\begin{crl}
\label{hfc for a reductively infinite operator 2}
Let $I$ be a set. Then every reductively transcendental bounded operator $A$ on $\t{C}_0(I,k)$ admits the holomorphic functional calculus in $\B_k(\t{C}_0(I,k))$ on $\n{A} \t{D}^1_k$.
\end{crl}

Thus every reductively infinite operator automatically admits a holomorphic functional calculus. Let $V$ be a strictly Cartesian Banach $k$-vector space and $A$ a bounded operator on $V$. Suppose $k = \C_p$ and $\sigma_{\B_{\C_p}(V)}(A) = O_{\C_p}$. As is referred in Example \ref{rigidity} (iii), we have $\C_p \ens{z} = \t{C}_{\t{rig}}(O_{\C_p},\C_p)$ by Theorem \ref{hfc for a reductively infinite operator}. Therefore in this case, the holomorphic functional calculus of $A$ in $\A$ on $\t{D}^1_k$ is the rigid continuous functional calculus. The following two examples might help a reader to understand the holomorphic functional calculus.

\begin{exm}
\label{shift}
Suppose that the valuation of $k$ is discrete and $\overline{k}$ is an infinite field. We consider the case $\A = \B_k(\t{C}_0(\Z,k))$. Let $A \in \A$ be the shift operator $(Af)(z) = f(z + 1)$. Then $A$ admits an inverse operator with norm $1$ given as the converse shift, and hence is unitary. Putting $M \coloneqq M_{\Z}(A)$, we have
\begin{eqnarray*}
  M =
  \left(
    \begin{array}{cccccc}
      \ddots & \vdots & \vdots & \vdots & \vdots & \t{\hspace{-.1em}\raisebox{-.7ex}{$\cdot$}\raisebox{.1ex}{$\cdot$}\raisebox{.9ex}{$\cdot$}} \\
      \cdots & 0 & 1 & 0 & 0 & \cdots \\
      \cdots & 0 & 0 & 1 & 0 & \cdots \\
      \cdots & 0 & 0 & 0 & 1 & \cdots \\
      \cdots & 0 & 0 & 0 & 0 & \cdots \\
      \t{\hspace{.1em}\raisebox{-.5ex}{$\cdot$}\raisebox{.2ex}{$\cdot$}\raisebox{.9ex}{$\cdot$}} & \vdots & \vdots & \vdots & \vdots & \ddots
    \end{array}
  \right).
\end{eqnarray*}
By the presentation of $M$, we obtain $\sigma_{\A}(A) = (O_k)^{\times}$. Since $\overline{k}$ is an infinite field, $A$ and $A^{-1}$ are reductively infinite. Applying Corollary \ref{hfc for a reductively infinite operator 2} to both $A$ and $A^{-1}$, we obtain the holomorphic functional calculus
\begin{eqnarray*}
  \iota_A \colon k \Ens{T,T^{-1}} & \to & \A \\
  F = \sum_{i \in \Z} F_i T^i & \mapsto & F(A) = \sum_{i \in \Z} F_i A^i
\end{eqnarray*}
of $A$ in $\A$ on the closed annulus of radius $1$. For any $F \in k \ens{T,T^{-1}}$, the matrix representation $F(M)$ of $F(A) \coloneqq \iota_A(F)$ is given as
\begin{eqnarray*}
  F(M) =
  \left(
    \begin{array}{cccccc}
      \ddots & \vdots & \vdots & \vdots & \vdots & \t{\hspace{-.1em}\raisebox{-.7ex}{$\cdot$}\raisebox{.1ex}{$\cdot$}\raisebox{.9ex}{$\cdot$}} \\
      \cdots & F_0 & F_1 & F_2 & F_3 & \cdots \\
      \cdots & F_{-1} & F_0 & F_1 & F_2 & \cdots \\
      \cdots & F_{-2} & F_{-1} & F_0 & F_1 & \cdots \\
      \cdots & F_{-3} & F_{-2} & F_{-1} & F_0 & \cdots \\
      \t{\hspace{.1em}\raisebox{-.5ex}{$\cdot$}\raisebox{.2ex}{$\cdot$}\raisebox{.9ex}{$\cdot$}} & \vdots & \vdots & \vdots & \vdots & \ddots
    \end{array}
  \right).
\end{eqnarray*}
This result is an analogue of the Fourier expansion of a function on $[0,1]$.
\end{exm}

\begin{exm}
\label{Toeplitz}
Suppose that the valuation of $k$ is discrete and $\overline{k}$ is an infinite field. We consider the case $\A = \B_k(\t{C}_0(\N,k))$. Let $A \in \A$ be the Toeplitz operator $(Af)(z) = f(z + 1)$. Then $A$ is an isometry with $\n{A} = 1$ and is not invertible because $A$ is not surjective. Putting $M \coloneqq M_{\N}(A)$, we have
\begin{eqnarray*}
  M =
  \left(
    \begin{array}{ccccc}
      0 & 1 & 0 & 0 & \cdots \\
      0 & 0 & 1 & 0 & \cdots \\
      0 & 0 & 0 & 1 & \cdots \\
      0 & 0 & 0 & 0 & \cdots \\
      \vdots & \vdots & \vdots & \vdots & \ddots
    \end{array}
  \right).
\end{eqnarray*}
By the presentation of $M$, we obtain $\sigma_{\A}(A) = O_k$. Since $\overline{k}$ is an infinite field, $A$ is reductively infinite. Applying Corollary \ref{hfc for a reductively infinite operator 2} to $A$, we obtain the holomorphic functional calculus
\begin{eqnarray*}
  \iota_A \colon k \Ens{T} & \to & \A \\
  F = \sum_{i \in \N} F_i T^i & \mapsto & F(A) = \sum_{i \in \N} F_i A^i
\end{eqnarray*}
of $A$ in $\A$ on the closed disc of radius $1$. For any $F \in k \ens{T}$, the matrix representation $F(M)$ of $F(A) \coloneqq \iota_A(F)$ is given as
\begin{eqnarray*}
  F(M) =
  \left(
    \begin{array}{ccccc}
      F_0 & F_1 & F_2 & F_3 & \cdots \\
      0 & F_0 & F_1 & F_2 & \cdots \\
      0 & 0 & F_0 & F_1 & \cdots \\
      0 & 0 & 0 & F_0 & \cdots \\
      \vdots & \vdots & \vdots & \vdots & \ddots
    \end{array}
  \right).
\end{eqnarray*}
The presentation of $F(M)$ for each $F \in k \ens{T}$ ensures that $\iota_A$ continuously extends to $O_k[[T]] \otimes_{O_k} k$, which is the $k$-algebra of sections on the open unit disc, with respect to the Fr\'echet topology on $O_k[[T]] \otimes_{O_k} k$ and the strong operator topology of $\A$.
\end{exm}

We note that the assumptions of $k$ in both examples above can be removed, because holomorphic functional calculi are compatible with restrictions of the scalar and base changes.

\subsection{Reductively Finite Operator}
\label{Reductively Finite Operator}

We deal with a reductively finite operator in this subsection. A reductively finite operator loses much information through the reduction. For example, we showed in Proposition \ref{reduction & operator reduction} that the operator reduction of a bounded operator of operator norm $1$ can be trivial, in which case the operator is trivially reductively finite. Even in the finite dimensional case, a matrix whose operator reduction is a diagonal matrix is not necessarily diagonalisable. Nevertheless, the reductive functional calculus of a reductively finite operator possesses information of a partition of the spectrum. Indeed, we will construct a decomposition of a reductively normal reductively finite operator associated to the eigenspace decomposition of the operator reduction in Theorem \ref{partition of unity}.

\begin{dfn}
Let $I$ be a set. A system $(P_i)_{i \in I} \in \A^I$ of idempotents in $\A$ is said to be a {\it partition of unity in $\A$} if $\set{i \in I}{P_i \neq 0}$ is a finite set, the essentially finite sum $\sum_{i \in I} P_i$ coincides with $1 \in \A$, and $\n{P_i}$ lies in $\ens{0,1}$ for any $i \in I$. A partition $(P_i)_{i \in I}$ of unity in $\A$ is said to be {\it reduced} if $\A \neq \ens{0}$ and $P_i \neq 0$ for any $i \in I$, is said to be {\it orthogonal} if $P_i P_j = 0$ for any $(i,j) \in I \times I$ with $i \neq j$, and is said to be {\it orthonormal} if it is reduced and orthogonal.
\end{dfn}

We note that in the case where the characteristic of $k$ is not $2$, a partition $(P_i)_{i \in I}$ of unity in $\A$ is orthogonal if and only if $P_i P_j = P_j P_i$ for any $(i,j) \in I \times I$.

\begin{prp}
\label{system of orthonormal projections}
Let $I$ be a set, and $(P_i)_{i \in I}$ an orthogonal partition of unity in $\A$. Then the equality $\n{A} = \max_{i \in I} \n{P_i A}$ holds for any $A \in \A$, and if $(P_i)_{i \in I}$ is reduced, then the equality $\n{\sum_{i \in I} c_i P_i} = \max_{i \in I} \v{c_i}$ holds for any $(c_i)_{i \in I} \in k^I$.
\end{prp}

\begin{proof}
The second assertion follows from the first assertion for $A = \sum_{i \in I} c_i P_i$, because the orthogonality of $(P_i)_{i \in I}$ ensures $P_{i_0} \sum_{i \in I} c_i P_i = c_{i_0} P_{i_0}$ for any $i_0 \in I$. Let $A \in \A$. The inequality $\n{A} \leq \max_{i \in I} \n{P_i A}$ follows from the equality $\sum_{i \in I} P_i A = (\sum_{i \in I} P_i) A = A$, and the inequality $\max_{i \in I} \n{P_i A} \leq \n{A}$ follows from the condition $\n{P_i} \in \ens{0,1}$ for each $i \in I$.
\end{proof}

\begin{prp}
\label{orthogonal property}
Let $V$ be a Banach $k$-vector space. A partition $(P_i)_{i \in I}$ of unity in $\B_k(V)$ induces a orthogonal decomposition $V = \bigoplus_{i \in I} P_i V = \bigoplus_{i \in I} \ker(1-P_i)$.
\end{prp}

\begin{proof}
Let $v \in V$. We have $\n{v} = \n{\sum_{i \in I} P_i v} \leq \max_{i \in I} \n{P_i v}$. On the other hand, we have $\max_{i \in I} \n{P_i v} \leq \max_{i \in I} \n{P_i} \ \n{v} \leq \n{v}$. Thus we obtain $\n{v} = \max_{i \in I} \n{P_i v}$.
\end{proof}

\begin{lmm}
\label{trivial idempotent}
An idempotent $P \in \A$ satisfies $\n{P} < 1$ if and only if $P = 0$.
\end{lmm}

\begin{proof}
The assertion follows from the inequality $\n{P} = \n{P^2} \leq \n{P}^2$.
\end{proof}

\begin{thm}
\label{partition of unity}
Let $A \in \A$ be a reductively finite element with $\n{A} = 1$. If $A$ admits the reductive functional calculus in $\A$, then there is a canonical orthogonal partition $(P_{A,\Lambda})_{\Lambda \in \overline{k}}$ of unity in $\A$ such that $P_{A,\Lambda}$ is non-zero for a $\Lambda \in \overline{k}$ if and only if $\Lambda \in \sigma_{\overline{\A}}(A + \A(1-))$, $P_{A,\Lambda}$ commutes with $A$ for any $\Lambda \in \overline{k}$, and the equalities
\begin{eqnarray*}
  \sigma_{\A}(P_{A,\Lambda} A) & = &
  \left\{
    \begin{array}{ll}
      \left( \Lambda \cap \sigma_{\A}(A) \right) \cup \ens{0} & \left( \sigma_{\overline{\A}}(A + \A(1-)) \neq \ens{\Lambda} \right) \\
      \sigma_{\A}(A) & \left( \sigma_{\overline{\A}}(A + \A(1-)) = \ens{\Lambda} \right)
    \end{array}
  \right. \\
  \sigma_{\A_{A,\Lambda}}(P_{A,\Lambda} A) & = &
  \left\{
    \begin{array}{ll}
      \Lambda \cap \sigma_{\A}(A) & \left( \Lambda \in \sigma_{\overline{\A}}(A + \A(1-)) \right) \\
      \emptyset & \left( \Lambda \in \overline{k} \backslash \sigma_{\overline{\A}}(A + \A(1-)) \right)
    \end{array}
  \right.
\end{eqnarray*}
hold for any $\Lambda \in \overline{k}$, where $\A_{A,\Lambda} \subset \A$ denotes the unital multiplicative $k$-vector subspace $P_{A,\Lambda}$ obtained as the closure of $P_{A,\Lambda} \A P_{A,\Lambda}$ for each $\Lambda \in \sigma_{\overline{\A}}(A + \A(1-))$.
\end{thm}

\begin{proof}
Put $\overline{A} \coloneqq A + \A(1-)$, $\sigma \coloneqq \sigma_{\A}(A)$, and $\overline{\sigma} \coloneqq \sigma_{\overline{\A}}(\overline{A})$. Set $P_{A,\Lambda} \coloneqq 0 \in \A$ for each $\Lambda \in \overline{k} \backslash \overline{\sigma}$. Put $n \coloneqq \# \overline{\sigma}$. Since $A$ is reductively finite, $n$ is a finite cardinal number. Since $A$ admits the reductive functional calculus, we have $n > 0$ by Remark \ref{non-empty spectrum}. If $n = 1$, then set $P_{A,\Lambda} \coloneqq 1 \in \A$ for a unique $\Lambda \in \overline{\sigma}$. In this case, $(P_{A,\Lambda})_{\Lambda \in \overline{k}}$ is a desired system. Suppose $n > 1$. Let $\Lambda \in \overline{\sigma}$. We denote by $\iota_{\overline{A}}$ the continuous functional calculus of $\overline{A}$ in $\overline{\A}$. The image of $\overline{k}{}^{\overline{\sigma}}$ by $\iota_{\overline{A}}$ is contained in $\Leb_{\overline{\A}}(\overline{A}) = \Leb(\overline{A})$ by Proposition \ref{continuous functions on a compact subset}, and $\Leb(\overline{A}) = \overline{k}[\overline{A}]$ is contained in the image of $\overline{\Leb(A)}$. Take a lift $P \in \Leb(A)(1)$ of $\iota_{\overline{A}}(\delta^{\overline{k}}_{\Lambda}) \in \overline{\A}$. The equality $(\delta^{\overline{k}}_{\Lambda})^2 = \delta^{\overline{k}}_{\Lambda}$ ensures $\n{P^2 - P} < 1$. Set $P_0 \coloneqq P$. We define a sequence $(P_j)_{j \in \N}$ on $\Leb(A)(1)$ by the recurrence relations $P_{j+1} = 3P_j^2 - 2P_j^3$ for each $j \in \N$. Then for any $j \in \N$, we have
\begin{eqnarray*}
  & & P_{j+1}^2 - P_{j+1} = (3P_j^2 - 2P_j^3)^2 - (3P_j^2 - 2P_j^3) = 4P_j^6 - 12P_j^5 + 9P_j^4 + 2P_j^3 - 3P_j^2 \\
  & = & (P_j^2 - P_j)(4P_j^4 -8P_j^3 + P_j^2 + 3P_j) = (P_j^2 - P_j)(P_j^2 - P_j)(4P_j^2 - 4P_j - 3) \\
  & = & -(P_j^2 - P_j)^2(3 - 4(P_j^2 - P_j))
\end{eqnarray*}
and hence $\n{P_{j+1}^2 - P_{j+1}} \leq \n{P_j^2 - P_j}^2$. Inductively on $j \in \N$, we obtain $\n{P_j^2 - P_j} \leq \n{P_0^2 - P_0}^{2^j} = \n{P^2 - P}^{2^j}$. It implies $\lim_{j \to \infty} \n{P_j^2 - P_j} = 0$. Moreover, we have $P_{j+1} - P_j = 3P_j^2 - 2P_j^3 -P_j = -(P_j^2 - P_j)(2P_j - 1)$, and hence $\n{P_{j+1} - P_j} \leq \n{P_j^2 - P_j} \ \n{2P_j - 1} \leq \n{P_j^2 - P_j}$. The completeness of $\A$ ensures that $(P_j)_{j \in \N}$ converges to some $P_{\Lambda} \in \A$. Since $\Leb(A)(1)$ is closed in $\A$, we have $P_{\Lambda} \in \Leb(A)(1)$. By the continuity of the norm, the multiplication, and the addition, we have $\n{P_{\Lambda}^2 - P_{\Lambda}} = \lim_{j \to \infty} \n{P_j^2 - P_j} = 0$, and hence $P_{\Lambda}$ is an idempotent. For any $j \in \N$, we have $\n{P_{j+1}} = \n{3P_j^2 - 2P_j^3} = \n{P_j + (P_j^2 - P_j)(-2P_j + 1)} = \n{P_j}$. Inductively on $j \in \N$, we obtain $\n{P_j} = \n{P_0} = \n{P} = 1$. It implies $\n{P_{\Lambda}} = \lim_{j \to \infty} \n{P_j} = 1$ by the continuity of the norm. We have $\n{P_{\Lambda} - P} = \n{\sum_{j = 0}^{\infty} (P_{j+1} - P_j)} \leq \sup_{j \in \N} \n{P_{j+1} - P_j} = \n{P_1 - P_0} = \n{(3P^2 - 2P^3) - P} = \n{-P(P - 1)(2P - 1)} \leq \n{P^2 - P} < 1$, and hence $P_{\Lambda} + \A(1-) = P + \A(1-) = \iota_{\A}(\delta^{(\overline{k})}_{\Lambda}) \neq 0$. Therefore we obtain $\n{P_{\Lambda}} = 1$.

\vspace{0.2in}
We verify that $P_{\Lambda}$ is independent of the choice of $P$. Let $Q_{\Lambda}$ be an idempotent in $\Leb(A)$ which lifts $\iota_{\A}(\delta^{(\overline{k})}_{\Lambda})$. Since $\Leb(A)$ is commutative, we have $(P_{\Lambda} - Q_{\Lambda})^4 = P_{\Lambda}^4 - 4P_{\Lambda}^3 Q_{\Lambda} + 6P_{\Lambda}^2 Q_{\Lambda}{}^2 - 4P_{\Lambda} Q_{\Lambda}^3 + Q_{\Lambda}^4 = P_{\Lambda} - 2P_{\Lambda} Q_{\Lambda} + Q_{\Lambda} = (P_{\Lambda} - Q_{\Lambda})^2$. It implies $\n{(P_{\Lambda} - Q_{\Lambda})^2} \leq \n{(P_{\Lambda} - Q_{\Lambda})}^2 < 1$ and $\n{(P_{\Lambda} - Q_{\Lambda})^2} = \n{(P_{\Lambda} - Q_{\Lambda})^4} \leq \n{(P_{\Lambda} - Q_{\Lambda})^2}^2$. Therefore we obtain $\n{(P_{\Lambda} - Q_{\Lambda})^2} = 0$. Moreover, we have
\begin{eqnarray*}
  & & (1 - P_{\Lambda})(1 - Q_{\Lambda}) = (1 - P_{\Lambda}) - P_{\Lambda} + P_{\Lambda} Q_{\Lambda} = (1 - P_{\Lambda}) - P_{\Lambda}^2 + P_{\Lambda} Q_{\Lambda} \\
  & = & (1 - P_{\Lambda}) + (Q_{\Lambda} - P_{\Lambda}) P_{\Lambda} = (1 - P_{\Lambda}) + (Q_{\Lambda} - P_{\Lambda}) \left( (3P_{\Lambda} - P_{\Lambda}) - (2P_{\Lambda} - P_{\Lambda}) \right) \\
  & = & (1 - P_{\Lambda}) + \left( 0 + 0 + 3(Q_{\Lambda} - P_{\Lambda}) P_{\Lambda} + P_{\Lambda} \right) - \left( 0 + 0 + 2(Q_{\Lambda} - P_{\Lambda})P_{\Lambda} + P_{\Lambda} \right) \\
  & = & (1 - P_{\Lambda}) + ((Q_{\Lambda} - P_{\Lambda}) + P_{\Lambda})^3 - ((Q_{\Lambda} - P_{\Lambda}) + P_{\Lambda})^2 \\
  & = & (1 - P_{\Lambda}) + (Q_{\Lambda}^3 - Q_{\Lambda}^2) = 1 - P_{\Lambda}
\end{eqnarray*}
and similarly $(1 - P_{\Lambda})(1 - Q_{\Lambda}) = 1 - Q_{\Lambda}$. We obtain $P_{\Lambda} = Q_{\Lambda}$. It implies that $P_{\Lambda}$ is a unique idempotent in $\Leb(A)$ which lifts $\iota_{\A}(\delta^{(\overline{k})}_{\Lambda})$, and hence is independent of the choice of $P$. Set $P_{A,\Lambda} \coloneqq P_{\Lambda}$. Since $\Leb(A)$ is commutative, $P_{A,\Lambda}$ commutes with $A$ by the construction for any $\Lambda \in \overline{k}$. The closure $\A_{A,\Lambda} \subset \A$ of $P_{A,\Lambda} \A P_{A,\Lambda}$ forms a Banach $k$-algebra with identity $P_{A,\Lambda}$. For any $B \in \A$ commuting with $P_{A,\Lambda}$, we have $P_{A,\Lambda} B = P_{A,\Lambda}^2 B = P_{A,\Lambda} B P_{A,\Lambda} \in \A_{A,\Lambda}$. In particular, we have $P_{A,\Lambda} A \in \A_{A,\Lambda}$.

\vspace{0.1in}
We have constructed a system $(P_{A,\Lambda})_{\Lambda \in \overline{k}}$ of idempotents in $\Leb(A)$ such that $P_{A,\Lambda}$ is non-zero for a $\Lambda \in \overline{k}$ if and only if $\Lambda \in \overline{\sigma}$. We verify that $(P_{A,\Lambda})_{\Lambda \in \overline{k}}$ forms an orthogonal partition of unity in $\A$. Let $(\Lambda, \Lambda') \in \overline{\sigma} \times \overline{\sigma}$ with $\Lambda \neq \Lambda'$. We show $P_{A,\Lambda} P_{A,\Lambda'} = 0$. We have $P_{A,\Lambda} P_{A,\Lambda'} + \A(1-) = (P_{A,\Lambda} + \A(1-))(P_{A,\Lambda'} + \A(1-)) = \iota_{\overline{A}}(\delta^{(\overline{k})}_{\Lambda}) \iota_{\overline{A}}(\delta^{(\overline{k})}_{\Lambda'}) = \iota_{\overline{A}}(\delta^{(\overline{k})}_{\Lambda} \delta^{(\overline{k})}_{\Lambda'}) = 0$ and hence $\n{P_{A,\Lambda} P_{A,\Lambda'}} < 1$. Moreover, $P_{A,\Lambda} P_{A,\Lambda'}$ is an idempotent, because $\Leb(A)$ is commutative. It implies $P_{A,\Lambda} P_{A,\Lambda'} = 0$ by Lemma \ref{trivial idempotent}.

\vspace{0.2in}
We have $(1 - \sum_{\Lambda \in \overline{k}} P_{A,\Lambda}) + \A(1-) = 1 - \sum_{\Lambda \in \overline{\sigma}} (P_{A,\Lambda} + \A(1-)) = \iota_{\overline{A}}(1) - \sum_{\Lambda \in \overline{\sigma}} \iota_{\overline{A}}(\delta^{(\overline{k})}_{\Lambda}) = \iota_{\overline{A}}(1 - \sum_{\Lambda \in \overline{\sigma}} \delta^{(\overline{k})}_{\Lambda}) = \iota_{\A}(0) = 0 \in \overline{\A}$, and hence $\n{1 - \sum_{\Lambda \in \overline{k}} P_{A,\Lambda}} < 1$. The equalities $P_{A,\Lambda} P_{A,\Lambda'} = 0$ for each $(\Lambda, \Lambda') \in \overline{k} \times \overline{k}$ with $\Lambda \neq \Lambda'$ ensure that $1 - \sum_{\Lambda \in \overline{k}} P_{A,\Lambda}$ is an idempotent, and hence coincides with $0$ by Lemma \ref{trivial idempotent}. Therefore $(P_{A,\Lambda})_{\Lambda \in \overline{k}}$ is an orthogonal partition of unity in $\A$.

\vspace{0.1in}
Let $\Lambda \in \overline{k}$. We verify the equalities on $\sigma_{\A}(P_{A,\Lambda} A)$ and $\sigma_{\A_{A,\Lambda}}(P_{A,\Lambda} A)$ in the assertion, which trivially hold if $\Lambda \notin \overline{\sigma}$. Suppose $\Lambda \in \overline{\sigma}$. We note that $P_{A,\Lambda} \in \A$ is not necessarily a central idempotent, and hence the $k$-linear homomorphism $\A \to \A_{A,\Lambda} \colon B \mapsto P_{A,\Lambda} B P_{A,\Lambda}$ is not necessarily a $k$-algebra homomorphism. Therefore the inclusion $\sigma_{\A_{A,\Lambda}}(P_{A,\Lambda} A) \subset \sigma_{\A}(P_{A,\Lambda} A)$ is not obvious. We have $\n{P_{A,\Lambda} A} \leq \n{A} = 1$ by Proposition \ref{system of orthonormal projections}, and hence $\sigma_{\A}(P_{A,\Lambda} A)$ and $\sigma_{\A_{A,\Lambda}}(P_{A,\Lambda} A)$ are contained in $O_k$ by Proposition \ref{spectrum -> bounded closed}. By the assumption $n > 1$, there is a $\Lambda' \in \overline{\sigma}$ with $\Lambda \neq \Lambda'$. The equality $P_{A,\Lambda} P_{A,\Lambda'} = 0$ with $P_{A,\Lambda'} \neq 0$ ensures $P_{A,\Lambda} \notin \A^{\times}$. Therefore we obtain $P_{A,\Lambda} A \notin \A^{\times}$ and $0 \in \sigma_{\A}(P_{A,\Lambda} A)$. We show $\sigma_{\A}(P_{A,\Lambda} A) \subset \sigma \cup \ens{0}$ and $\sigma_{\A_{A,\Lambda}}(P_{A,\Lambda} A) \subset \sigma$. Let $\lambda \in O_k \backslash \sigma$. Then $A - \lambda$ is invertible in $\A$ by definition. We have $(P_{A,\Lambda} A - \lambda P_{A,\Lambda})(P_{A,\Lambda} (A - \lambda)^{-1} P_{A,\Lambda}) = P_{A,\Lambda}^3 = P_{A,\Lambda} \in \A_{A,\Lambda}$ and $(P_{A,\Lambda} (A - \lambda)^{-1} P_{A,\Lambda})(P_{A,\Lambda} A - \lambda P_{A,\Lambda}) = P_{A,\Lambda}^3 = P_{A,\Lambda} \in \A_{A,\Lambda}$ because $P_{A,\Lambda}$ commutes with $A$. It implies $P_{A,\Lambda} A - \lambda P_{A,\Lambda} \in \A_{A,\Lambda}^{\times}$ and $\lambda \in O_k \backslash \sigma_{\A_{A,\Lambda}}(P_{A,\Lambda} A)$. In addition if $\lambda \neq 0$, then we have
\begin{eqnarray*}
  & & (P_{A,\Lambda} A - \lambda) \left( P_{A,\Lambda} (A - \lambda)^{-1} P_{A,\Lambda} + (1 - P_{A,\Lambda}) \lambda^{-1} \right) \\
  & = & (P_{A,\Lambda} (A - \lambda) + (1 - P_{A,\Lambda}) \lambda) \left( P_{A,\Lambda} (A - \lambda)^{-1} P_{A,\Lambda} + (1 - P_{A,\Lambda}) \lambda^{-1} \right) \\
  & = & P_{A,\Lambda}^3 + P_{A,\Lambda} (1 - P_{A,\Lambda}) (A - \lambda) \lambda^{-1} + (1 - P_{A,\Lambda}) P_{A,\Lambda} \lambda (A - \lambda)^{-1} P_{A,\Lambda} A + (1 - P_{A,\Lambda})^2 \\
  & = & P_{A,\Lambda} + 0 + 0 + (1 - P_{A,\Lambda}) = 1
\end{eqnarray*}
and 
\begin{eqnarray*}
  & & \left( P_{A,\Lambda} (A - \lambda)^{-1} P_{A,\Lambda} + (1 - P_{A,\Lambda}) \lambda^{-1} \right)(P_{A,\Lambda} A - \lambda) \\
  & = & \left( P_{A,\Lambda} (A - \lambda)^{-1} P_{A,\Lambda} + (1 - P_{A,\Lambda}) \lambda^{-1} \right)(P_{A,\Lambda} (A - \lambda) + (1 - P_{A,\Lambda}) \lambda) \\
  & = & P_{A,\Lambda}^3 + P_{A,\Lambda} (A - \lambda)^{-1} P_{A,\Lambda} (1 - P_{A,\Lambda}) \lambda + (1 - P_{A,\Lambda}) P_{A,\Lambda} \lambda^{-1} (A - \lambda) + (1 - P_{A,\Lambda})^2 \\
  & = & P_{A,\Lambda} + 0 + 0 + (1 - P_{A,\Lambda}) = 1.
\end{eqnarray*}
It implies $P_{A,\Lambda} A - \lambda \in \A^{\times}$ and $\lambda \in O_k \backslash \sigma_{\A}(P_{A,\Lambda} A)$. Therefore we obtain $\sigma_{\A}(P_{A,\Lambda} A) \subset \sigma \cup \ens{0}$ and $\sigma_{\A_{A,\Lambda}}(P_{A,\Lambda} A) \subset \sigma$.

\vspace{0.2in}
We show $\sigma_{\A}(P_{A,\Lambda} A) \subset \Lambda$ and $\sigma_{\A_{A,\Lambda}}(P_{A,\Lambda} A) \subset \Lambda$. Let $\lambda \in O_k \backslash \Lambda$. Since $\lambda + m_k$ does not coincide with $\Lambda$, we have $(P_{A,\Lambda} A - \lambda P_{A,\Lambda}) + \A(1-) = \iota_{\overline{A}}(\delta^{\overline{k}}_{\Lambda} z^{(\overline{k})}_{\overline{\sigma}} - (\lambda + \A(1-)) \delta^{\overline{k}}_{\Lambda}) = \iota_{\overline{A}}(\Lambda \delta^{\overline{k}}_{\Lambda} - (\lambda + \A(1-)) \delta^{\overline{k}}_{\Lambda}) = (\Lambda - (\lambda + \A(1-))) (P_{A,\Lambda} + \A(1-)) \in \overline{\A_{A,\Lambda}}{}^{\times}$, and hence $P_{A,\Lambda} A - \lambda P_{A,\Lambda} \in \A_{A,\Lambda}(1)^{\times}$ by Lemma \ref{inverse and reduction}. It implies $\lambda \in O_k \backslash \sigma_{\A_{A,\Lambda}}(P_{A,\Lambda} A)$. Let $B \in \A_{A,\Lambda}^{\times}$ be the inverse of $P_{A,\Lambda} A - P_{A,\Lambda} \lambda$ in $\A_{A,\Lambda}$. The equalities $P_{A,\Lambda} (1 - P_{A,\Lambda}) = 0$ and $(1 - P_{A,\Lambda}) P_{A,\Lambda} = 0$ ensure that $B_0 (1 - P_{A,\Lambda}) = (1 - P_{A,\Lambda}) B_0 = 0$ for any $B_0 \in P_{A,\Lambda} \A P_{A,\Lambda}$. Since $B$ lies in $\A_{A,\Lambda}$, $B$ is the limit of some sequence on $P_{A,\Lambda} \A P_{A,\Lambda}$, and the continuity of the multiplication on $\A$ ensures $B(1 - P_{A,\Lambda}) = (1 - P_{A,\Lambda})B = 0$. In addition if $\lambda \neq 0$, then we have $(P_{A,\Lambda} A - \lambda)(B - \lambda^{-1}(1 - B)) = (B - \lambda^{-1}(1 - B))(P_{A,\Lambda} A - \lambda) = 1$. It implies $P_{A,\Lambda} A - \lambda \in \A^{\times}$ and $\lambda \in O_k \backslash \sigma_{\A}(P_{A,\Lambda} A)$. Therefore we obtain $\sigma_{\A}(P_{A,\Lambda} A) \subset \Lambda \cup \ens{0}$ and $\sigma_{\A_{A,\Lambda}}(P_{A,\Lambda} A) \subset \Lambda$.

\vspace{0.2in}
We have verified $0 \in \sigma_{\A}(P_{A,\Lambda} A)$, $\sigma_{\A}(P_{A,\Lambda} A) \subset (\Lambda \cap \sigma) \cup \ens{0}$, and $\sigma_{\A_{A,\Lambda}}(P_{A,\Lambda} A) \subset \Lambda \cap \sigma$. Therefore it suffices to verify $\Lambda \cap \sigma \subset \sigma_{\A}(P_{A,\Lambda} A)$ and $\Lambda \cap \sigma \subset \sigma_{\A_{A,\Lambda}}(P_{A,\Lambda} A)$. Let $\lambda \in \Lambda \cap \sigma$. First, we show $\lambda \in \sigma_{\A}(P_{A,\Lambda} A)$. We have $P_{A,\Lambda} A - \lambda$ is invertible in $\A$ by definition. Assume $\lambda \notin \sigma_{\A}(P_{A,\Lambda} A)$. We show that $P_{A,\Lambda'} A - \lambda$ is invertible in $\A$ for any $\Lambda' \in \overline{\sigma}$. Let $\Lambda' \in \overline{\sigma}$. If $\Lambda' = \Lambda$, then we have $\lambda \notin \sigma_{\A}(P_{A,\Lambda'} A)$, and hence $P_{A,\Lambda'} A - \lambda$ is invertible in $\A$. Suppose $\Lambda' \neq \Lambda$. Since $0$ lies in $\sigma_{\A}(P_{A,\Lambda} A)$, we have $\lambda \neq 0$. Since $\lambda$ lies in $\Lambda \neq \Lambda'$, we obtain $\lambda \notin (\Lambda' \cap \sigma) \cup \ens{0}$, and hence $\lambda \notin \sigma_{\A}(P_{A,\Lambda'} A)$. Therefore $P_{A,\Lambda'} A - \lambda$ is invertible in $\A$ for any $\Lambda' \in \overline{\sigma}$. The equalities
\begin{eqnarray*}
  & & (A - \lambda) \left( \sum_{\Lambda' \in \overline{\sigma}} P_{A,\Lambda'} (P_{A,\Lambda'} A - \lambda)^{-1} P_{A,\Lambda'} \right) \\
  & = & \left( \sum_{\Lambda' \in \overline{\sigma}} P_{A,\Lambda'} (A - \lambda) \right) \left( \sum_{\Lambda' \in \overline{\sigma}} P_{A,\Lambda'} (P_{A,\Lambda'} A - \lambda)^{-1} P_{A,\Lambda'} \right) \\
  & = & \left( \sum_{\Lambda' \in \overline{\sigma}} P_{A,\Lambda'} (P_{A,\Lambda'} A - \lambda) \right) \left( \sum_{\Lambda' \in \overline{\sigma}} P_{A,\Lambda'} (P_{A,\Lambda'} A - \lambda)^{-1} P_{A,\Lambda'} \right) \\
  & = & \sum_{\Lambda' \in \overline{\sigma}} \sum_{\Lambda'' \in \overline{\sigma}} P_{A,\Lambda'} P_{A,\Lambda''} (P_{A,\Lambda'}A - \lambda) (P_{A,\Lambda''} A - \lambda)^{-1} P_{A,\Lambda'} = \sum_{\Lambda' \in \overline{\sigma}} P_{A,\Lambda'}^2 = \sum_{\Lambda' \in \overline{\sigma}} P_{A,\Lambda'} = 1
\end{eqnarray*}
and
\begin{eqnarray*}
  & & \left( \sum_{\Lambda' \in \overline{\sigma}} P_{A,\Lambda'} (P_{A,\Lambda'} A - \lambda)^{-1} P_{A,\Lambda'} \right)(A - \lambda) \\
  & = & \left( \sum_{\Lambda' \in \overline{\sigma}} P_{A,\Lambda'} (P_{A,\Lambda'} A - \lambda)^{-1} P_{A,\Lambda'} \right) \left( \sum_{\Lambda' \in \overline{\sigma}} P_{A,\Lambda'} (A - \lambda) \right) \\
  & = & \left( \sum_{\Lambda' \in \overline{\sigma}} P_{A,\Lambda'} (P_{A,\Lambda'} A - \lambda)^{-1} P_{A,\Lambda'} \right) \left( \sum_{\Lambda' \in \overline{\sigma}} P_{A,\Lambda'} (P_{A,\Lambda'} A - \lambda) \right) \\
  & = & \sum_{\Lambda' \in \overline{\sigma}} \sum_{\Lambda'' \in \overline{\sigma}} P_{A,\Lambda'} (P_{A,\Lambda'} A - \lambda)^{-1} P_{A,\Lambda'} P_{A,\Lambda''} (P_{A,\Lambda''} A - \lambda) = \sum_{\Lambda' \in \overline{\sigma}} P_{A,\Lambda'}^2 = \sum_{\Lambda' \in \overline{\sigma}} P_{A,\Lambda'} = 1
\end{eqnarray*}
imply $A - \lambda \in \A^{\times}$, and it contradicts $\lambda \in \Lambda \cap \sigma \subset \sigma$. As a consequence, we obtain $\lambda \in \sigma_{\A}(P_{A,\Lambda} A)$.

\vspace{0.2in}
Secondly, we show $\lambda \in \sigma_{\A_{A,\Lambda}}(P_{A,\Lambda} A)$. Assume $\lambda \notin \sigma_{\A_{A,\Lambda}}(P_{A,\Lambda} A)$. We prove that $P_{A,\Lambda'} A - \lambda P_{A,\Lambda'}$ is invertible in $\A_{A,\Lambda}$ for any $\Lambda' \in \overline{\sigma}$. Let $\Lambda' \in \overline{\sigma}$. If $\Lambda' = \Lambda$, then we have $\lambda \notin \sigma_{\A}(P_{A,\Lambda'} A)$, and hence $P_{A,\Lambda'} A - \lambda P_{A,\Lambda'}$ is invertible in $\A_{A,\Lambda}$. Suppose $\Lambda' \neq \Lambda$. Since $\lambda$ lies in $\Lambda \neq \Lambda'$, we obtain $\lambda \notin \Lambda' \cap \sigma$, and hence $\lambda \notin \sigma_{\A_{A,\Lambda}}(P_{A,\Lambda'} A)$. Therefore $P_{A,\Lambda'} A - \lambda P_{A,\Lambda'}$ is invertible in $\A_{A,\Lambda'}$ for any $\Lambda' \in \overline{\sigma}$. Let $\Lambda' \in \overline{\sigma}$. We denote by $B_{\Lambda'}$ the inverse of $P_{A,\Lambda'} A - \lambda P_{A,\Lambda'}$ in $\A_{A,\Lambda'}$, i.e.\ a unique element in $\A_{A,\Lambda'}$ satisfying $(P_{A,\Lambda'} A - \lambda P_{A,\Lambda'}) B_{\Lambda'} = B_{\Lambda'} (P_{A,\Lambda'} A - \lambda P_{A,\Lambda'}) = P_{A,\Lambda'}$. Since $P_{A,\Lambda'}$ is the identity of $\A_{A,\Lambda'}$, we have $P_{A,\Lambda'}B_{\Lambda'}P_{A,\Lambda'} = B_{\Lambda'}$. The equalities
\begin{eqnarray*}
  & & (A - \lambda) \left( \sum_{\Lambda' \in \overline{\sigma}} B_{\Lambda'} \right) = \left( \sum_{\Lambda' \in \overline{\sigma}} P_{A,\Lambda'} (A - \lambda) \right) \left( \sum_{\Lambda' \in \overline{\sigma}} P_{A,\Lambda'} B_{\Lambda'} P_{A,\Lambda'} \right) \\
  & = & \sum_{\Lambda' \in \overline{\sigma}} \sum_{\Lambda'' \in \overline{\sigma}} P_{A,\Lambda'} P_{A,\Lambda''} (A - \lambda) B_{\Lambda''} P_{A,\Lambda''} = \sum_{\Lambda' \in \overline{\sigma}} P_{A,\Lambda'}^2 (A - \lambda) B_{\Lambda'} P_{A,\Lambda'} \\
  & = & \sum_{\Lambda' \in \overline{\sigma}} (P_{A,\Lambda'} A - \lambda P_{A,\Lambda'}) B_{\Lambda'} P_{A,\Lambda'} = \sum_{\Lambda' \in \overline{\sigma}} P_{A,\Lambda'}^2 = \sum_{\Lambda' \in \overline{\sigma}} P_{A,\Lambda'} = 1
\end{eqnarray*}
and
\begin{eqnarray*}
  & & \left( \sum_{\Lambda' \in \overline{\sigma}} B_{\Lambda'} \right) (A - \lambda) = \left( \sum_{\Lambda' \in \overline{\sigma}} P_{A,\Lambda'} B_{\Lambda'} P_{A,\Lambda'} \right) \left( \sum_{\Lambda' \in \overline{\sigma}} P_{A,\Lambda'} (A - \lambda) \right) \\
  & = & \sum_{\Lambda' \in \overline{\sigma}} \sum_{\Lambda'' \in \overline{\sigma}} P_{A,\Lambda'} B_{\Lambda'} P_{A,\Lambda'} P_{A,\Lambda''} (A - \lambda) = \sum_{\Lambda' \in \overline{\sigma}} P_{A,\Lambda'} B_{\Lambda'} P_{A,\Lambda'}^2 (A - \lambda) \\
  & = & \sum_{\Lambda' \in \overline{\sigma}} P_{A,\Lambda'} B_{\Lambda'} (P_{A,\Lambda''} A - \lambda P_{A,\Lambda'}) = \sum_{\Lambda' \in \overline{\sigma}} P_{A,\Lambda'}^2 = \sum_{\Lambda' \in \overline{\sigma}} P_{A,\Lambda'} = 1
\end{eqnarray*}
imply $A - \lambda \in \A^{\times}$, and it contradicts $\lambda \in \Lambda \cap \sigma \subset \sigma$. As a consequence, we obtain $\lambda \in \sigma_{\A_{A,\Lambda}}(P_{A,\Lambda} A)$. We conclude $\sigma_{\A}(P_{A,\Lambda} A) = (\Lambda \cap \sigma) \cup \ens{0}$ and $\sigma_{\A_{A,\Lambda}}(P_{A,\Lambda} A) =
\Lambda \cap \sigma$.
\end{proof}

Let $A \in \A$ be a reductively finite element with $\n{A} = 1$ admitting the reductive functional calculus in $\A$. As is shown in Example \ref{spectrum of affinoids}, $\sigma_{\overline{\A}}(A + \A(1-))$ depends not only on $A$ but also on $\A$. Although the definition of $(P_{A,\Lambda})_{\Lambda \in \overline{k}}$ uses $\sigma_{\overline{\A}}(A + \A(1-))$, it is independent of $\A$. Indeed, the condition that $A$ admits the reductive functional calculus in $\A$ ensures the stability of $\sigma_{\overline{\A}}(A + \A(1-))$ under changes of $\A$ containing $A$. We note that $(P_{A,\Lambda})_{\Lambda \in \overline{k}}$ is obtained as a system in $\Leb(A) \subset \A$, but $\Leb(A)$ does not necessarily satisfy the condition that $A$ admits the reductive functional calculus in $\Leb(A)$.

\begin{prp}
\label{rn + ln <-> n}
Let $A \in \A$ be a reductively finite element with $\n{A} = 1$ admitting the reductive functional calculus in $\A$. Then $A$ admits the rigid continuous functional calculus in $\A$ if and only if $P_{A,\Lambda} A$ admits the rigid continuous functional calculus in $\A_{A,\Lambda}$ for any $\Lambda \in \overline{k}$.
\end{prp}

\begin{proof}
Suppose that $A$ admits the rigid continuous functional calculus $\iota_A$ in $\A$. Since $(P_{A,\Lambda})_{\Lambda \in \overline{k}}$ is a system of idempotents in $\Leb(A)$, it is contained in the image of $\Leb(z_{\sigma_{\A}(A)}) \subset \t{C}_{\t{rig}}(\sigma_{\A}(A),k)$ by $\iota_A$. For any $\Lambda \in \overline{k} \backslash \sigma_{\overline{\A}}(A + \A(1-))$, $P_{A,\Lambda} A$ admits the rigid continuous functional calculus in $\A_{A,\Lambda} = \ens{0}$ because of the equalities $P_{A,\Lambda} = 0$, $\sigma_{\A_{A,\Lambda}}(P_{A,\Lambda} A) = \sigma_{\ens{0}}(0) = \emptyset$, and $\t{C}(\emptyset,k) = \ens{0}$. Let $(E_{\Lambda})_{\Lambda \in \overline{k}}$ denote the partition of unity in $\t{C}_{\t{rig}}(\sigma_{\A}(A),k)$ given as the preimage of $(P_{A,\Lambda})_{\Lambda \in \overline{k}}$ by $\iota_{A}$. For each $\Lambda \in \overline{k}$, we denote by $\t{supp}(E_{\Lambda}) \subset \sigma_{\A}(A)$ the measurable subset $\set{\lambda \in k}{E_{\lambda}(\lambda) = 1}$, and by $\iota^{(\Lambda)}_! \colon \t{C}_{\t{rig}}(\t{supp}(E_{\Lambda}),k) \hookrightarrow \t{C}_{\t{rig}}(\sigma_{\A}(A),k)$ the zero-extension in Theorem \ref{functoriality}. We show $\t{supp}(E_{\Lambda}) = \Lambda \cap \sigma_{\A}(A)$ for any $\Lambda \in \overline{k}$. Let $\Lambda \in \overline{k}$. Since $\iota_A \circ \iota^{(\Lambda)}_!$ is an isometric multiplicative $k$-linear homomorphism sending $E_{\Lambda}$ to $P_{A,\Lambda}$, it induces an isometric $k$-algebra homomorphism $\iota_{P_{A,\Lambda} A} \colon \t{C}_{\t{rig}}(\t{supp}(E_{\Lambda}),k) \hookrightarrow \A_{A,\Lambda}$ sending $z^{(k)}_{\t{supp}(E_{\Lambda})}$ to $\iota_A(E_{\Lambda} z^{(k)}_{\sigma_{\A}(A)}) = P_{A,\Lambda} A$. Therefore we obtain $\Lambda \cap \sigma_{\A}(A) = \sigma_{\A_{A,\Lambda}}(P_{A,\Lambda} A) \subset \sigma_{\t{C}_{\t{rig}}(\t{supp}(E_{\Lambda}),k)}(z^{(k)}_{\t{supp}(E_{\Lambda})}) = \t{supp}(E_{\Lambda})$ by Theorem \ref{invertible rigid continuous function} and Proposition \ref{functoriality 2} (i). Let $\lambda \in \t{supp}(E_{\Lambda})$. Since $\t{supp}(E_{\Lambda})$ is a subset of $\sigma_{\A}(A)$, $A - \lambda$ is not invertible in $\A$. Since $(P_{A,\Lambda'})_{\Lambda' \in \overline{k}}$ is a partition of unity in $\A$ consisting of idempotents commuting with $A$ and each other, $P_{A,\Lambda_0}(A - \lambda)$ is not invertible in $\A_{A,\Lambda_0}$ for some $\Lambda_0 \in \overline{k}$. It implies $\lambda \in \sigma_{\A_{A,\Lambda_0}}(P_{A,\Lambda_0} A) = \Lambda_0 \cap \sigma_{\A}(A) \subset \t{supp}(E_{\Lambda_0})$. Since $(P_{A,\Lambda'})_{\Lambda' \in \overline{k}}$ is a partition of unity in $\A$, so is $(E_{\Lambda'})_{\Lambda' \in \overline{k}}$ in $\t{C}_{\t{rig}}(\sigma_{\A}(A),k)$. Therefore we obtain $\sigma_{\A}(A) = \bigsqcup_{\Lambda' \in \overline{k}} \t{supp}(E_{\Lambda'})$. Since $\t{supp}(E_{\Lambda})$ shares a common element $\lambda$ with $\t{supp}(E_{\Lambda_0})$, we get $\Lambda = \Lambda_0$, and hence $\lambda \in \Lambda_0 \cap \sigma_{\A}(A) = \Lambda \cap \sigma_{\A}(A)$. It implies $\t{supp}(E_{\Lambda}) = \Lambda \cap \sigma_{\A}(A)$. Since $E_{\Lambda}$ is an idempotent, we obtain $E_{\Lambda} = \delta^{(k)}_{\Lambda \cap \sigma_{\A}(A),\sigma_{\A}(A)}$. Therefore $\iota_{P_{A,\Lambda}}$ is the rigid continuous functional calculus of $P_{A,\Lambda} A$ in $\A_{A,\Lambda}$.

\vspace{0.2in}
Suppose that $P_{A,\Lambda} A$ admits the rigid continuous functional calculus $\iota_{P_{A,\Lambda} A}$ in $\A_{A,\Lambda}$ for any $\Lambda \in \overline{k}$. Let $\varphi \colon \prod_{\Lambda \in \overline{k}} \A_{A,\Lambda} \to \A$ denote the $k$-linear homomorphism given by the addition. Then $\varphi$ sends $(P_{A,\Lambda} A)_{\Lambda \in \overline{k}}$ and $(P_{A,\Lambda} A)_{\Lambda \in \overline{k}}$ to $1$ and $A$ respectively, and is an isometry by Proposition \ref{system of orthonormal projections}, because $(P_{A,\Lambda})_{\Lambda \in \overline{k}}$ is a partition of unity in $\A$ consisting of idempotents commuting with $A$. Since $(P_{A,\Lambda})_{\Lambda \in \overline{k}}$ is orthogonal, $\varphi$ is a multiplicative. Therefore $\varphi$ is an isometric $k$-algebra homomorphism. Let $\iota_A \colon \t{C}_{\t{rig}}(\sigma_{\A}(A),k) \to \A$ denote the composite of the isometric $k$-algebra isomorphism $\t{C}_{\t{rig}}(\sigma_{\A}(A),k) \to \prod_{\Lambda \in \overline{k}} \t{C}_{\t{rig}}(\Lambda \cap \sigma_{\A}(A),k)$ induced by the restriction maps in Theorem \ref{functoriality}, the isometric $k$-algebra homomorphism $\prod_{\Lambda \in \overline{k}} \t{C}_{\t{rig}}(\Lambda \cap \sigma_{\A}(A),k) \to \prod_{\Lambda \in \overline{k}} \A_{A,\Lambda}$ given as the direct product of $(\iota_{P_{A,\Lambda} A})_{\Lambda \in \overline{k}}$, and $\varphi$. Then $\iota_A$ is an isometric $k$-algebra homomorphism sending $z^{(k)}_{\sigma_{\A}(A)} = \sum_{\Lambda \in \overline{k}} \iota^{(\Lambda)}_!(z^{(k)}_{\Lambda \cap \sigma_{\A}(A)})$ to $\sum_{\Lambda \in \overline{k}} \iota_{P_{A,\Lambda} A}(z^{(k)}_{\Lambda \cap \sigma_{\A}(A)}) = \sum_{\Lambda \in \overline{k}} P_{A,\Lambda} A = A$, and is the rigid continuous functional calculus of $A$ in $\A$.
\end{proof}

Henceforth, we assume that the valuation of $k$ is discrete and $\n{\A}$ is contained in $\v{k}$. We fix a uniformiser $\varpi_k \in k^{\times}$. For an $A \in \A \backslash \ens{0}$, we put $v(A) \coloneqq \log_{\v{\varpi_k}} \n{A} \in \Z$. We call the map $v(\cdot) \colon \A \to \Z$ {\it the additive norm on $\A$}, which is independent of the choice of $\varpi_k$. We set $\Nor_{-1}(\A) \coloneqq \A$, define maps $P_0 \colon \Nor_{-1}(\A) \to \A$ and $r_0 \colon \Nor_{-1}(\A) \to \A$ by setting $P_0(A) \coloneqq 1$ and $r_0(A) = P_0(A) A = A$ for each $A \in \A$, and put $\A_0^{(A)} \coloneqq \A = P_0(A) \A P_0(A)$ for each $A \in \Nor_{-1}(\A)$. An $A \in \A$ is said to be {\it naive of level $0$ at $0$ in $\A$} if one of the following conditions holds:
\begin{itemize}
\item[(i)] The equality $A = 0$ holds.
\item[(ii)] The reduction $O_k \twoheadrightarrow \overline{k}$ induces a surjective map $\sigma_{\A}(c^{-1}A) \twoheadrightarrow \sigma_{\overline{\A}}(c^{-1}A + \A(1-))$ for a $c \in k^{\times}$ with $\n{A} = \v{c}$, and $A$ admits the reductive continuous functional calculus in $\A$.
\end{itemize}
We denote by $\Nor_0(\A) \subset \A$ the subset of naive elements of level $0$ at $0$ in $\A$, define maps $P(\cdot,\A) = P_1 \colon \Nor_0(\A) \to \A$ and $r_1 \colon \Nor_0(\A) \to \A$ by setting $P(0,\A) = P_1(0) \coloneqq 0$, $P(A,\A) = P_1(A) \coloneqq P_{\varpi_k^{-v(A)} A, 0}$ for each $A \in \Nor_0(\A) \backslash \ens{0}$, and $r_1(A) \coloneqq P_1(A) A = P(A,\A) A$ for each $A \in \Nor_0(\A)$, and denote by $\A_1^{(A)}$ the closure of the unital multiplicative $k$-linear subspace $P_1(A) \A P_1(A) \subset \A$ with the identity $P_1(A)$ for each $A \in \Nor_0(\A)$. Inductively on $i \in \N$, we set $\Nor_{i+1}(\A) \coloneqq r_1^{-1}(\Nor_i(\A)) \subset \Nor_i(\A)$, define maps $P_{i+2} \colon \Nor_{i+1}(\A) \to P_{i+1}(A) \A$ and $r_{i+2} \colon \Nor_{i+1}(\A) \to \A$ by setting $P_{i+2}(A) \coloneqq P(P_{i+1}(A) A, \A_{i+1}^{(A)}) P_{i+1}(A)$ and $r_{i+2} \coloneqq r_{i + 1}(r_1(A))$, and denote by $\A_{i+2}^{(A)}$ the closure of the multiplicative $k$-linear subspace $P_{i+2}(A) \A P_{i+2}(A) \subset \A$ with the identity $P_{i+2}(A)$ for each $A \in \Nor_{i+1}(\A)$. Then we have $r_{i+1}(A) = P_{i+1}(A) A$ and $P_{i+2}(A) P_{i+1}(A) = P_{i+1}(A) P_{i+2}(A) = P_{i+2}(A)$ for any $(i,A) \in \N \times \Nor_i(\A)$. Let $\lambda \in \sigma_{\A}(A)$. An $A \in \A$ is said to be {\it naive of level $i$ at $\lambda$ in $\A$} for an $i \in \N$ if $A - \lambda \in \Nor_i(\A)$, and is said to be {\it naive at $\lambda$ in $\A$} if $A$ is naive of level $i$ at $\lambda$ in $\A$ for any $i \in \N$. An $A \in \A$ is said to be {\it naive in $\A$} if $A$ is naive at $\lambda$ in $\A$ for any $\lambda \in \sigma_{\A}(A)$.

\begin{prp}
\label{naivety}
For any non-empty compact subset $\sigma \subset k$, the equality $\n{\t{C}(\sigma,k)} = \v{k}$ holds, and $z_{\sigma}$ is naive in $\t{C}(\sigma,k)$.
\end{prp}

\begin{proof}
The first implication follows from the maximal modulus principle of a continuous function. Let $(i,\lambda) \in \N \times \sigma$. We verify that $z_{\sigma}$ is naive of level $i$ at $\lambda \in \sigma$. We have $\sigma_{\t{C}(\sigma,k)}(z_{\sigma} - \lambda) = \sigma - \lambda \coloneqq \set{\lambda' - \lambda}{\lambda' \in \sigma}$, and the parallel transform $\sigma \to \sigma - \lambda \colon \lambda' \mapsto \lambda' - \lambda$ induces an isometric $k$-algebra isomorphism $\t{C}(\sigma,k) \to \t{C}(\sigma - \lambda, k)$ sending $z_{\sigma} - \lambda$ to $z_{\sigma - \lambda}$. Therefore replacing $\sigma$ by $\sigma - \lambda$, we may assume $\lambda = 0 \in \sigma$. If $z_{\sigma} \equiv 0$, i.e.\ $\sigma = \ens{0}$, then $z_{\sigma}$ is naive at $0$ by definition. Otherwise, replacing $\sigma$ by $c^{-1} \sigma$ for a $c \in k^{\times}$ with $\n{z_{\sigma}} = \v{c}$, we may assume $\n{z_{\sigma}} = 1$. We denote by $\overline{\sigma} \subset \overline{k}$ the image of $\sigma$ by the reduction $O_k \twoheadrightarrow \overline{k}$. We have $\sigma_{\overline{\t{C}(\sigma,k)}}(z_{\sigma} + \t{C}(\sigma,k)(1-)) = \overline{\sigma}$ by Corollary \ref{spectrum and reduction 3}. Theorem \ref{n implies rn 2} ensures that $z_{\sigma}$ is naive of level $0$ at $0$ in $\t{C}(\sigma,k)$, and $(P_{z_{\sigma},\Lambda})_{\Lambda \in \overline{k}}$ is the system of idempotents corresponding to the partition of $\sigma$ associated to the canonical projection $\sigma \twoheadrightarrow \overline{\sigma}$. Under the identification between $\t{C}(\sigma,k)_{z_{\sigma},0} \subset \t{C}(\sigma,k)$ and $\t{C}(m_k \cap \sigma, k)$ given by the restriction map, $r_1(z_{\sigma}) = P_{z_{\sigma},0} z_{\sigma} \in \t{C}(\sigma,k)$ corresponds to $z_{m_k \cap \sigma}$. Therefore the induction on $i$ works after the replacement of $\sigma$ by $m_k \cap \sigma$. We conclude that $z_{\sigma}$ is naive at $0$ in $\t{C}(\sigma,k)$.
\end{proof}

\begin{lmm}
\label{isolation}
Let $A \in \A$ be an element with $\n{A} = 1$ admitting the reductive functional calculus in $\A$. Then $P_{A,0} A$ coincides with $0$ if and only if $m_k \cap \sigma_{\A}(A)$ is contained in $\ens{0}$.
\end{lmm}

\begin{proof}
The necessary implication is obvious by the definition of $P_{A,0}$. If $P_{A,0} A = 0$, then we obtain $m_k \cap \sigma_{\A}(A) \subset \sigma_{\A}(P_{A,0} A) = \ens{0}$ after noting that $m_k$ is $0 \in \overline{k}$ as a coset.
\end{proof}

\begin{lmm}
\label{decreasing the radius}
The inequality $\n{r_{i+1}(A)} \leq \v{\varpi_k}^{i+1} \n{A}$ holds for any $i \in \N$ and $A \in \Nor_i(\A)$.
\end{lmm}

\begin{proof}
It suffices to verify the inequality in the case $i = 0$ and $\n{A} = 1$ by the construction of $r_{i+1}$. In this case, we have $A \neq 0$ and hence $A$ admits the reductive continuous functional calculus $\iota_{A + \A(1-)}$ in $\A$. If $0$ lies in $\overline{k} \backslash \sigma_{\overline{\A}}(A + \A(1-))$, then we have $\n{r_1(A)} = \n{P_{A,0} A} = \n{0} = 0 < \v{\varpi_k}$. Suppose $0 \in \sigma_{\overline{\A}}(A + \A(1-))$. By the construction of $P_{A,0}$, we have $P_{A,0} A + \A(1-) = (P_{A,0} + \A(1-))(A + \A(1-)) = \iota_{A + \A(1-)}(\delta^{(\overline{k})}_0) \iota_{A + \A(1-)}(z_{\sigma_{\overline{\A}}(A + \A(1-))}) = \iota_{A + \A(1-)}(\delta^{(\overline{k})}_0 z_{\sigma_{\overline{\A}}(A + \A(1-))}) = \iota_{A + \A(1-)}(0) = 0$, and hence $\n{r_1(A)} = \n{P_{A,0} A} < 1$. It implies $\n{r_1(A)} \leq \v{\varpi_k}$.
\end{proof}

\begin{lmm}
\label{approximation}
Let $A \in \A$ be an element with $\n{A} = 1$ of naive of level $0$ at every $\lambda \in \A(1)$. Then $A - \sum_{\Lambda \in \overline{k}} P_1(A - \lbrack \Lambda \rbrack)$ coincides with $\sum_{\Lambda \in \overline{k}} P_{A,\Lambda} \lbrack \Lambda \rbrack$, and satisfies $\n{A - \sum_{\Lambda \in \overline{k}} P_1(A - \lbrack \Lambda \rbrack)} \leq \v{\varpi_k}$, where $\lbrack \cdot \rbrack \colon \overline{k} \hookrightarrow k$ denote the Teichm\"uller embedding.
\end{lmm}

\begin{proof}
We have $P_1(A - \lbrack \Lambda \rbrack) = P_{A - \lbrack \Lambda \rbrack,0}(A - \lbrack \Lambda \rbrack) = P_{A,\Lambda} (A - \lbrack \Lambda \rbrack)$ for any $\Lambda \in \overline{k}$, and hence $A - \sum_{\Lambda \in \overline{k}} P_1(A - \lbrack \Lambda \rbrack) = A - (\sum_{\Lambda \in \overline{k}} P_{A,\Lambda} (A - \lbrack \Lambda \rbrack)) = \sum_{\Lambda \in \overline{k}} \lbrack \Lambda \rbrack P_{A,\Lambda}$. In order to verify the second assertion, it suffices to show $\n{P_1(A')} \leq \v{\varpi_k}$ for any $A' \in \Nor_0(\A)$ with $\n{A'} = 1$. Since $P_1(A') = P_{A',0} A'$ is a lift of $\iota_{A' + \A(1-)}(\delta^{(\overline{k})}_0)(A' + \A(1-)) = \iota_{A' + \A(1-)}(\delta^{(\overline{k})}_0 z^{(\overline{k})}_{\sigma_{\A}(A')}) = 0$, it lies in $\A(1-)$. By the assumption $\n{A} \subset \v{k}$, we conclude $\n{P_1(A')} \leq \v{\varpi_k}$.
\end{proof}

\begin{thm}
\label{reduction and functional calculus}
Suppose that $k$ is a local field. An $A \in \A$ admits the continuous functional calculus in $\A$ if and only if $A$ is naive in $\A$.
\end{thm}

\begin{proof}
We remark that every bounded closed subset $\sigma \subset k$ is compact and every $A \in \A \backslash \ens{0}$ is reductively finite, because $k$ is a local field. Let $A \in \A$. Suppose that $A$ admits the continuous functional calculus in $\A$. Then $A$ is naive by Corollary \ref{spectrum and reduction 4}, Theorem \ref{n implies rn 2}, and Proposition \ref{rn + ln <-> n}.

\vspace{0.2in}
Suppose that $A$ is naive in $\A$. Let $c \colon \sigma_{\A}(A) \to k$ be a locally constant function. We construct a substitution $c(A) \in \A$. Since $\sigma_{\A}(A)$ is compact, there is an $r > 0$ such that $c$ is constant on each closed disc of radius $r$ in $\sigma_{\A}(A)$. Replacing $r$ by a smaller one, we may assume $r \in \v{k^{\times}}$, because the valuation of a local field is non-trivial. Let $\lambda \in \sigma_{\A}(A)$. Since $A$ is naive at $\lambda$, there is an $i \in \N$ such that $\n{r_i(A - \lambda)} \leq r$ by Lemma \ref{decreasing the radius}. We denote by $i_{\lambda} \in \N$ the smallest integer among such $i$'s. Lemma \ref{isolation} ensures that $\n{r_{i_{\lambda}}(A - \lambda)}$ coincides with $0$ only when $\lambda$ is isolated in $\sigma_{\A}(A)$ by the construction of $r_{i_{\lambda}}$, and hence the closed disc $\set{\lambda' \in \sigma_{\A}(A)}{\v{\lambda' - \lambda} \leq \n{r_{i_{\lambda}}(A - \lambda)}}$ is an open subset of $\sigma_{\A}(A)$. We obtain an open covering $\U \coloneqq \set{\set{\lambda' \in \sigma_{\A}(A)}{\v{\lambda' - \lambda} \leq \n{r_{i_{\lambda}}(A - \lambda)}}}{\lambda \in \sigma_{\A}(A)}$ of $\sigma_{\A}(A)$. Since $\sigma_{\A}(A)$ is compact, $\U$ admits a finite subcovering $\U_0 \subset \U$. For each $U \in \U_0$, fix a centre $\lambda_U \in U$. We define $c(A)$ as $\sum_{U \in \U_0} c(\lambda_U)P_{i_{\lambda_U}}(A - \lambda_U) \in \A$. By Proposition \ref{system of orthonormal projections}, we have $\n{c(A)} = \max{}_{U \in \U_0} \v{c(\lambda_U)} = \n{c}$. The sum of a partition of unity is $1$, and hence $c(A)$ is stable under a refinement of $\U_0$ by the compatibility of the spectra and the reductions in the definition of the naivety. It implies that $c(A)$ is independent of the choice of $r > 0$ and $\U_0$. If $c$ is a constant map, then $c(A)$ coincides with $c \in k$ by definition.

\vspace{0.2in}
Let $c,c' \colon \sigma_{\A}(A) \to k$ be two locally constant functions. Take finite disjoint open coverings of $\sigma_{\A}(A)$ by closed discs in the definition of $c(A)$ and $c'(A)$. Refining them, we may assume that they are a common covering. Then the equalities $(c + c')(A) = c(A) + c'(A)$ and $(cc')(A) = c(A)c'(A)$ hold by the construction. Therefore the substitution $c \mapsto c(A)$ gives a $k$-algebra homomorphism $\iota_A$ from the $k$-subalgebra $\t{c}(\sigma_{\A}(A),k) \subset \t{C}(\sigma_{\A}(A),k)$ of locally constant functions to $\A$, which is an isometry with respect to the restriction of the supremum norm of $\t{C}(\sigma_{\A}(A),k)$. Since $\sigma_{\A}(A)$ is a totally disconnected compact Hausdorff topological space, $\t{c}(\sigma_{\A}(A),k)$ is dense in $\t{C}(\sigma_{\A}(A),k)$, and hence $\iota_A$ can be uniquely extended to be an isometry $\iota_A \colon \t{C}(\sigma_{\A}(A),k) \to \A$. It sends a sequence on $\t{c}(\sigma_{\A}(A),k)$ converging to $z_{\sigma_{\A}(A)}$ in $\t{C}(\sigma_{\A}(A),k)$ to a sequence converging to $A$ in $\A$, because Lemma \ref{approximation} ensures that every naive element $A' \in \A$ can be approximated by the linear combination of idempotents in modulo $\varpi_k \A(\n{A'})$. We conclude that $\iota_A$ is the continuous functional calculus of $A$ in $\A$.
\end{proof}

\begin{crl}
\label{reduction and functional calculus 2}
Suppose that $k$ is a local field. Let $V$ be a strictly Cartesian Banach $k$-vector space, and $A$ a bounded operator on $V$. Then $A$ is normal if and only if $A$ is naive in $\B_k(V)$.
\end{crl}

Let $V$ be a strictly Cartesian Banach $k$-vector space, and $A$ a bounded operator on $V$. Computation of the naivety of $A$ includes recursive calculi with infinite steps unless $V$ is of finite dimension. Therefore practical use of the criterion for the normality using the naivety needs an assumption that $A$ admits a good structure like a ``fractal''. See the following examples.

\begin{exm}
Suppose that $k$ is a finite extension of $\Q_p$. Set $\A \coloneqq \B_k(\t{C}_0(\Z,k))$. We consider the elements $A_1,A_2 \in \A$ whose matrix representations with respect to the canonical orthonormal Schauder basis of $\t{C}_0(\Z,k)$ are given as
\begin{eqnarray*}
  & & M_1 =
  \left(
    \begin{array}{ccccccc}
      \ddots & \vdots & \vdots & \vdots & \vdots & \vdots & \t{\hspace{-.1em}\raisebox{-.7ex}{$\cdot$}\raisebox{.1ex}{$\cdot$}\raisebox{.9ex}{$\cdot$}} \\
      \cdots & -2 & p & p^2 & p^3 & p^4 & \cdots \\
      \cdots & 0 & -1 & p & p^2 & p^3 & \cdots \\
      \cdots & 0 & 0 & 0 & p & p^2 & \cdots \\
      \cdots & 0 & 0 & 0 & 1 & p & \cdots \\
      \cdots & 0 & 0 & 0 & 0 & 2 & \cdots \\
      \t{\hspace{.1em}\raisebox{-.5ex}{$\cdot$}\raisebox{.2ex}{$\cdot$}\raisebox{.9ex}{$\cdot$}} & \vdots & \vdots & \vdots & \vdots & \vdots & \ddots
    \end{array}
  \right), \\
  & & M_2 =
  \left(
    \begin{array}{ccccccc}
      \ddots & \vdots & \vdots & \vdots & \vdots & \vdots & \t{\hspace{-.1em}\raisebox{-.7ex}{$\cdot$}\raisebox{.1ex}{$\cdot$}\raisebox{.9ex}{$\cdot$}} \\
      \cdots & -2 & p & 0 & 0 & 0 & \cdots \\
      \cdots & 0 & -1 & p & 0 & 0 & \cdots \\
      \cdots & 0 & 0 & 0 & p & 0 & \cdots \\
      \cdots & 0 & 0 & 0 & 1 & p & \cdots \\
      \cdots & 0 & 0 & 0 & 0 & 2 & \cdots \\
      \t{\hspace{.1em}\raisebox{-.5ex}{$\cdot$}\raisebox{.2ex}{$\cdot$}\raisebox{.9ex}{$\cdot$}} & \vdots & \vdots & \vdots & \vdots & \vdots & \ddots
    \end{array}
  \right)
\end{eqnarray*}
respectively. Let $h \in \ens{1,2}$. We have $\n{A_h} = 1$, $\sigma_{\A}(A_h) = \Z_p$, and $\sigma_{\overline{\A}}(A_h + \A(1-)) = \F_p \subset \overline{k}$. The reduction of $A_h$ is naturally identified with the operator reduction of $M_h$ through the isomorphism in Proposition \ref{reduction & operator reduction}, and admits the continuous functional calculus because its matrix representation is a diagonal matrix. Moreover, $P_{A_h, i + m_k} A_h \in \A_{A_h, j + m_k}$ acting on $P_{A_h, i + m_k} \t{C}_0(\Z,k) = \t{C}_0(i + p \Z, k)$ admits the matrix representation quite similar to that of $A_h \in \A$ acting on $\t{C}_0(\Z,k)$ for each $i \in \N \cap [0,p-1]$. A simple repetition of computation ensures that all the operators appearing in the recursive process in the verification of the naivety of $A$ have a similar matrix representation. Therefore $A_h$ is naive, and hence admits the continuous functional calculus $\iota_{A_h} \colon \t{C}(\Z_p,k) \to \A$. For any continuous map $F \colon \Z_p \to k$, the matrix representations $F(M_1)$ and $F(M_2)$ of $\iota_{A_1}(F) = F(A_1)$ and $\iota_{A_2}(F) = F(A_2)$ are given as
\begin{eqnarray*}
  & & F(M_1) =
  \left(
    \begin{array}{ccccccc}
      \ddots & \vdots & \vdots & \vdots & \vdots & \vdots & \t{\hspace{-.1em}\raisebox{-.7ex}{$\cdot$}\raisebox{.1ex}{$\cdot$}\raisebox{.9ex}{$\cdot$}} \\
      \cdots & F(-2) & \partial F(-2)p & \partial F(-2)p^2 & \partial F(-2)p^3 & \partial F(-2)p^4 & \cdots \\
      \cdots & 0 & F(-1) & \partial F(-1)p & \partial F(-1)p^2 & \partial F(-1)p^3 & \cdots \\
      \cdots & 0 & 0 & F(0) & \partial F(0)p & \partial F(0)p^2 & \cdots \\
      \cdots & 0 & 0 & 0 & F(1) & \partial F(1)p & \cdots \\
      \cdots & 0 & 0 & 0 & 0 & F(2) & \cdots \\
      \t{\hspace{.1em}\raisebox{-.5ex}{$\cdot$}\raisebox{.2ex}{$\cdot$}\raisebox{.9ex}{$\cdot$}} & \vdots & \vdots & \vdots & \vdots & \vdots & \ddots
    \end{array}
  \right), \\
  & & F(M_2) =
  \left(
    \begin{array}{ccccccc}
      \ddots & \vdots & \vdots & \vdots & \vdots & \vdots & \t{\hspace{-.1em}\raisebox{-.7ex}{$\cdot$}\raisebox{.1ex}{$\cdot$}\raisebox{.9ex}{$\cdot$}} \\
      \cdots & F(-2) & \partial F(-2)p & \frac{\partial^2 F(-2)}{2!}p^2 & \frac{\partial^3 F(-2)}{3!}p^3 & \frac{\partial^4 F(-2)}{4!}p^4 & \cdots \\
      \cdots & 0 & F(-1) & \partial F(-1)p & \frac{\partial^2 F(-1)}{2!}p^2 & \frac{\partial^3 F(-1)}{3!}p^3 & \cdots \\
      \cdots & 0 & 0 & F(0) & \partial F(0)p & \frac{\partial^2 F(0)}{2!}p^2 & \cdots \\
      \cdots & 0 & 0 & 0 & F(1) & \partial F(1)p & \cdots \\
      \cdots & 0 & 0 & 0 & 0 & F(2) & \cdots \\
      \t{\hspace{.1em}\raisebox{-.5ex}{$\cdot$}\raisebox{.2ex}{$\cdot$}\raisebox{.9ex}{$\cdot$}} & \vdots & \vdots & \vdots & \vdots & \vdots & \ddots
    \end{array}
  \right)
\end{eqnarray*}
respectively, where we set $\partial G(n) \coloneqq G(n) - G(n - 1)$ for each $n \in \Z_p$ and $G \in \t{C}(\Z_p,k)$. This result is obviously analogous to the holomorphic functional calculus of the shift operator in Example \ref{shift}. We remark that the difference $(n!)^{-1} \partial^n F$ appears in continuous functional calculi, while the differential coefficient $F_n = (n!)^{-1} d^nF/dT^n$ does in holomorphic functional calculi. The differential operator $\frac{d}{dz}$ is continuous in the class of {\it rigid} analytic functions $F \in k \ens{z}$, and the difference operator $\partial$ is continuous in the class of continuous functions $F \in \t{C}(\Z_p,k)$, which are {\it naive} by Proposition \ref{naivety}.
\end{exm}

\subsection{Reduction Algorithm for Matrices}
\label{Reduction Algorithm for Matrices}

Continuing from Theorem \ref{reduction and functional calculus}, we assume that $k$ is a local field. Following the philosophy of the repetitive reduction derived from Theorem \ref{reduction and functional calculus}, we construct an explicit algorithm for a criterion of the diagonalisability of a matrix by a unitary matrix. Let $V$ be a strictly Cartesian Banach $k$-vector space, and $A$ a bounded operator on $V$ with $\n{A} = 1$. We verified in Proposition \ref{rn <-> rfc} that $A$ is reductively normal if and only if $A$ admits the reductive functional calculus in $\overline{\B_k(V)}$. If $\sigma_{\overline{\B_k(V)}}(A + \B_k(V)(1-))$ contains at least two elements, then Proposition \ref{rn + ln <-> n} helps us to determine whether $A$ is normal or not. If $\sigma_{\overline{\B_k(V)}}(A + \B_k(V)(1-))$ consists of a single element $\Lambda_0 \in \overline{k}$, then we have $P_{A,\Lambda_0} = 1$ and $P_{A,\Lambda} = 0$ for any $\Lambda \in \overline{k} \backslash \ens{\Lambda_0}$. In this case, Proposition \ref{rn + ln <-> n} gives no information about the normality of $A$. In order to apply Proposition \ref{rn + ln <-> n} to a construction of a criterion of the normality of $A$, we need to replace $A$ by a suitable operator possessing the information of the normality of $A$ so that the spectrum of the operator reduction contains at least two elements in the following way.

\begin{prp}
\label{reductively scalar}
Let $A \in \A$ be an element with $\n{A} = 1$ admitting the reductive functional calculus in $\A$. Then a $c \in O_k$ satisfies $\sigma_{\overline{\A}}(A + \A(1-)) = \ens{c + m_k}$ if and only if $\n{A - c} < 1$, in which case $A$ lies in $\A(1)^{\times}$ and $c$ lies in $O_k^{\times}$.
\end{prp}

\begin{proof}
First, suppose $\sigma_{\overline{\A}}(A + \A(1-)) = \ens{c + m_k}$. Since $A$ admits the reductive functional calculus in $\A$, we have $A + \A(1-) = z^{(\overline{k})}_{\ens{c + m_k}} = c + m_k$. It implies $\n{A - c} < 1$ and $\v{c} = 1$. Next, suppose $\n{A - a} < 1$. Then we obtain $(A + \A(1-)) - (c + m_k) = (A - c) + A(1-) = 0 \in \overline{\A}$, and hence $\ens{c + m_k} = \sigma_{\overline{\A}}(c + m_k) = \sigma_{\overline{\A}}(A + \A(1-))$ after noting $\overline{\A} \neq \ens{0}$. The inequality $\n{A - c} < 1$ ensures $\v{c} = 1$. It implies $0 \notin \sigma_{\overline{\A}}(A + \A(1-))$, and hence $A$ lies in $\A(1)^{\times}$ by Lemma \ref{inverse and reduction}.
\end{proof}

\begin{crl}
\label{reductively scalar 2}
Let $V$ be a strictly Cartesian Banach $k$-vector space, and $A$ a bounded operator on $V$. Then $\sigma_{\overline{\B_k(V)}}(A + \B_k(V)(1-))$ consists of a single element if and only if $A$ is reductively scalar.
\end{crl}

\begin{proof}
The assertion follows from Proposition \ref{reduction & operator reduction} and Proposition \ref{reductively scalar}.
\end{proof}

\begin{lmm}
\label{reductively non-scalar}
For any $M = (M_{i,j})_{i,j=1}^n \in \t{M}_n(k)$, $M - M_{1,1}$ is reductively scalar if and only if $M$ is scalar.
\end{lmm}

\begin{proof}
Since the $(1,1)$-entry of $M - M_{1,1}$ is $0$, the assertion follows from the fact that a matrix which possesses $0$ in a diagonal entry is reductively scalar if and only if it is the zero-matrix.
\end{proof}

\begin{thm}
\label{unitarily diagonalisable}
Let $M = (M_{i,j})_{i,j=1}^n \in \t{M}_n(k)$ be a non-scalar matrix. Then $M$ is diagonalisable by a unitary matrix in $k$ if and only if $\Pi_{k^n}(M - M_{1,1}) \in \t{M}_n(\overline{k})$ is diagonalisable in $\overline{k}$ and the matrix representation of the restriction of $M - M_{1,1}$ to the $M$-stable subspace $\ker \prod_{\lambda \in \Lambda} (M - M_{1,1} - \lambda)^n \subset k^n$ with respect to an orthonormal basis is diagonalisable by a unitary matrix in $k$ for any $\Lambda \in \Sigma(M - M_{1,1})$, in which case $\sigma_{\t{M}_n(\overline{k})}(\Pi_{k^n}(M - M_{1,1}))$ contains at least two elements
\end{thm}

\begin{proof}
We note that by \cite{BGR84} 2.5.1/5, every subspace of $k^n$ is a strictly Cartesian Banach $k$-vector space and admits an orthonormal basis. The first assertion follows from Proposition \ref{reduction & operator reduction}, Proposition \ref{n <-> ud}, and Proposition \ref{rn + ln <-> n}, and the second assertion follows from Corollary \ref{reductively scalar 2} and Lemma \ref{reductively non-scalar}.
\end{proof}

As a consequence, Theorem \ref{unitarily diagonalisable} yields an algorithm for a criterion of the diagonalisability of a matrix by a unitary matrix, because the dimension of the representation space strictly decreases in each recursive step of the decomposition of the representation space induced by the eigenspace decomposition of the operator reduction.

\vspace{0.4in}
\addcontentsline{toc}{section}{Acknowledgements}
\noindent {\Large \bf Acknowledgements}
\vspace{0.1in}

I am profoundly grateful to Takeshi Tsuji for giving advices in seminars. I appreciate helpful personal discussions with Anatoly N.\ Kochubei. I am thankful to my colleagues for discussions. I express my gratitude to family for their deep affection.

\end{document}